\documentclass[10pt]{article}

\usepackage{amsmath}
\usepackage{amssymb}
\usepackage{latexsym}
\usepackage{amsthm}

\newcommand{\init}{\big\vert_{t = 0}}
 \newcommand{\half}{\frac{1}{2}}
\newcommand{\abs}[1]{\left\vert #1 \right\vert}
\newcommand{\bigabs}[1]{\bigl\vert #1 \bigr\vert}

\newcommand{\norm}[1]{\left\Vert #1 \right\Vert}
\newcommand{\bignorm}[1]{\bigl\Vert #1 \bigr\Vert}

\newcommand{\spacetimenorm}[3]{\norm{#1}_{H^{#2,#3}}}

\newcommand{\Spacetimenorm}[3]{\norm{#1}_{\scrH^{#2,#3}}}

\newcommand{\Sobnorm}[2]{\norm{#1}_{H^{#2}}}
\newcommand{\bigSobnorm}[2]{\bignorm{#1}_{H^{#2}}}
\newcommand{\Sobdotnorm}[2]{\norm{#1}_{\dot{H}^{#2}}}
\newcommand{\bigSobdotnorm}[2]{\bignorm{#1}_{\dot{H}^{#2}}}
\newcommand{\Lpnorm}[2]{\norm{#1}_{L^{#2}}}
\newcommand{\bigLpnorm}[2]{\bignorm{#1}_{L^{#2}}}

\newcommand{\twonorm}[2]{\norm{#1}_{L^2#2}}
\newcommand{\bigtwonorm}[2]{\bignorm{#1}_{L^2#2}}

\newcommand{\mixednorm}[3]{\norm{#1}_{L_{t}^{#2}(L_{x}^{#3})}}

\newcommand{\Mixednorm}[3]{\norm{#1}_{\LL_{t}^{#2}(\LL_{x}^{#3})}}
\newcommand{\Mixed}[2]{\LL_{t}^{#1}(\LL_{x}^{#2})}
\newcommand{\bigMixednorm}[3]{\bignorm{#1}_{\LL_{t}^{#2}(\LL_{x}^{#3})}}
\newcommand{\bigmixednorm}[3]{\bignorm{#1}_{L_{t}^{#2}(L_{x}^{#3})}}

\newcommand{\C}{\mathbb{C}}

\newcommand{\scrH}{\mathcal{H}}
\newcommand{\LL}{\mathcal{L}}

\newcommand{\Proj}{\mathcal{P}}
\newcommand{\R}{\mathbb{R}}
\newcommand{\X}{\mathcal{X}}
\newcommand{\Y}{\mathcal{Y}}
\newcommand{\Z}{\mathcal{Z}}

\newcommand{\Schwartz}{\mathcal{S}}
\newcommand{\Fourier}{\mathcal{F}}

\newcommand{\hypwt}[2]{\bigabs{ \abs{#1} - \abs{#2} }}

\DeclareMathOperator{\diag}{diag}
\DeclareMathOperator{\curl}{curl}

\newtheorem{theorem}{Theorem}
\newtheorem{proposition}{Proposition}
\newtheorem{lemma}{Lemma}

\newtheorem*{plaintheorem}{Theorem}

\theoremstyle{definition}

\theoremstyle{remark}
\newtheorem{remark}{Remark}
\newtheorem*{plainremark}{Remark}
\newtheorem*{plainremarks}{Remarks}

\title{ALMOST OPTIMAL LOCAL WELL-POSEDNESS OF THE MAXWELL-KLEIN-GORDON 
EQUATIONS IN $1+4$ DIMENSIONS}

\author{Sigmund Selberg\footnote{Current address: Inst.\@ f.\@
Mathematik, Univ.\@ Wien, Strudlhofgasse 4, A-1090 Wien, Austria}\\
Department of Mathematics\\ Johns Hopkins University
\\Baltimore, MD 21218}

\date{}

\begin{document}

\maketitle

\begin{abstract}
We prove that the Maxwell-Klein-Gordon system on $\R^{1+4}$
relative to the Coulomb gauge is
locally well-posed for initial data in $H^{1+\varepsilon}$ for all
$\varepsilon > 0$.
This builds on previous work by Klainerman and Machedon \cite{Kl-Ma5}
who proved the corresponding result, with the additional restriction 
of small-norm data, for a model problem obtained
by ignoring the elliptic features of the system, as well as cubic terms.
\end{abstract}
\section{Introduction}
The purpose of this paper is to prove local well-posedness (LWP)
of the Maxwell-Klein-Gordon (MKG) equations on
$\R^{1+4}$, relative to the Coulomb gauge, for initial data in $H^{1+\varepsilon}$,
any $\varepsilon > 0$. This result is optimal in the sense that the
critical Sobolev exponent for MKG on $\R^{1+4}$ is $s_c = 1$,
and one does not expect well-posedness in $H^{s}$ for $s$
below this critical value; see the introduction in \cite{Kl-Ta}
and section \ref{Scaling} below,
where we also make some remarks on the open question of well-posedness
in the critical data norm $H^{1}$.

The analogous result for a hyperbolic model problem,
obtained from the MKG system \eqref{MKG} below
by setting the non-dynamical variable $A_0 \equiv 0$
and ignoring all cubic terms, was proved by Klainerman-Machedon
\cite{Kl-Ma5}, for small-norm initial data. That result was reproved,
using different norms, and without any smallness assumption
on the data, in the recent survey article \cite{Kl-Se}.
The proof given there also used some ideas from
\cite{Kl-Ta}, where the corresponding model problem for the
Yang-Mills equation is considered.

The present work builds on the treatment of the
model problem in \cite{Kl-Se}: To obtain \emph{a priori}
estimates on solutions of MKG with the requisite regularity,
we complement the bilinear estimates proved there with estimates
for cubic terms, and terms involving the non-dynamical variable,
which satisfies an elliptic equation. It should be emphasized that 
the difficulty is to obtain LWP when $s$ is very close to $s_{c} = 1$.
If $s$ is sufficiently large, one can prove LWP by much simpler methods 
than those employed here. See section \ref{Scaling} and Remark 
\ref{SpacesRemark} in section \ref{Notation}. 

Our method here can be modified\footnote{We do not prove this 
here, but hope to address it in a separate paper 
dealing with the Yang-Mills equations on $\R^{1+4}$ in Coulomb
as well as temporal gauge. Note that Yang-Mills essentially contains MKG as 
a special case.} to treat the full Yang-Mills 
system in $\R^{1+4}$, proving LWP in 
$H^{1+\varepsilon}$, but only for initial data with small norm.
This extends the result of Klainerman-Tataru \cite{Kl-Ta} on
a model equation for Yang-Mills.
The reason for the small-norm restriction is that the elliptic equation in the 
Yang-Mills system relative to the Coulomb gauge is far more 
complicated than the one for MKG, and not in general globally solvable.
To avoid this problem one can include the elliptic variable in the Picard
iteration. Then to close the iteration one must assume small-norm 
data, since there is no way of compensating for large data by letting 
the existence time go to zero, as one can for an iteration involving only 
hyperbolic equations in a subcritical regime.
Of course, using Picard iteration for an elliptic equation seems 
somewhat contrived. A better approach for Yang-Mills on $\R^{1+4}$ may be
to work in the temporal gauge, as Tao \cite{TaoYM} has successfully done for
the case of $\R^{1+3}$. We hope to address this in a future paper.

Most of the previous work on MKG has been in dimension $1+3$. Let us 
summarize the known results for this case.
LWP in the energy norm $H^{1}$ was proved by Klainerman 
and Machedon \cite{Kl-Ma0.2}. By conservation of the MKG energy, their
result implies global well-posedness. In particular, they recovered an
earlier global regularity result of Eardley and Moncrief \cite{E-M}
for smooth data.
Cuccagna \cite{Cu} proved LWP for small-norm data in $H^{s}$,
$s > 3/4$. For $1+3$ dimensions,
the critical regularity is $s_{c} = 1/2$,
but the question of LWP below $s = 3/4$ remains open.
In both \cite{Kl-Ma0.2} and \cite{Cu} the Coulomb gauge is used.
More recently, Tao \cite{TaoYM} has proved small-norm LWP 
for $s > 3/4$ using the temporal gauge, for the more general 
Yang-Mills equations.

Our method here can be used to remove the small-norm restriction in the
result of Cuccagna. The essential reason for this limitation in \cite{Cu}
is that the elliptic variable was included in the iteration. If instead one solves
the elliptic equation and reduces to a purely hyperbolic system, as 
we do here, this obstruction is removed, and one can get a large data LWP
result. A crucial fact needed to make this work is that in the 
Klainerman-Machedon bilinear estimates used by Cuccagna, the 
space-time derivative $\abs{D_{t,x}}^{-a}$ acting on the product
can be replaced by 
$\abs{D_{x}}^{-a}$, as observed in \cite{Se2} (cf. also the remark in the 
Appendix), rendering unnecessary the decomposition in Fourier space 
used in \cite{Cu}. See also Remark \ref{CuRemark} in section 
\ref{EllipticEstimates}.

Finally, we remark that our proof should generalize without difficulity
to the higher dimensional case of MKG on $\R^{1+n}$ with $n \ge 5$,
giving LWP for $s > s_{c} = \tfrac{n-2}{2}$.
In fact, the difficulty of the problem decreases with increasing dimension.
\subsection{The Maxwell-Klein-Gordon system}
The Klein-Gordon equation can be derived as a 
relativistic analogue of the Schr\"odinger equation for a free particle. 
It is obtained from the relativistic energy-momentum relation
$E^{2} = \mathbf p^{2} c^{2} + m^{2} c^{4}$, where $E$ is the energy 
of the particle, $m > 0$ its rest mass, $\mathbf p$ its momentum and $c$
the light speed. Setting $c = 1$ from now on, and applying the quantum
mechanical principle of replacing classical quantities by operators:
\begin{itemize}
  \item[] Energy $\qquad E \longrightarrow i \frac{\partial}{\partial 
  t}$,
  \item[] Momentum $\qquad \mathbf p \longrightarrow \frac{1}{i} \nabla$,
\end{itemize}
one obtains the free Klein-Gordon equation
\begin{equation}\label{FreeKG}
  \square \phi = m^{2} \phi,
\end{equation}
where $\phi(t,x) \in \C$ and $\square = \partial_\mu \partial^\mu
= - \partial_{t}^{2} + \Delta$ is the wave operator
on $\R^{1+n}$. Here we use relativistic coordinates $t = x^{0}, 
x^{1}, \dots, x^{n}$ on the Minkowski spacetime $\R^{1+n}$ with the metric
$\diag(-1,1,\dots,1)$; indices are raised and lowered relative to this 
metric, and the Einstein summation convention is in effect: roman indices
$j,k,\dots$ run from $1$ to $n$, greek indices $\mu,\nu,\dots$ from
$0$ to $n$. We write $\partial_{\mu}$ for $\tfrac{\partial}{\partial 
x^{\mu}}$, and $\partial_{t} = \partial_{0}$. We shall use
$\Re z$ and $\Im z$ to denote the real and imaginary parts of $z \in \C$.

The coupling of \eqref{FreeKG} to an 
electromagnetic field represented by a potential $A_{\mu}(t,x) \in \R$ is achieved 
by the so-called minimal substitution
$$
  \partial_{\mu} \longrightarrow D_{\mu} = \partial_{\mu} + i A_{\mu},
$$
where $i A_{\mu}$ acts as a multiplication operator. This gives
\begin{equation}\label{KG}
  D_{\mu} D^{\mu} \phi = m^{2} \phi.
\end{equation}
which is the Klein-Gordon equation.
It has an associated current density
\begin{equation}\label{KGcurrent}
  j_{\mu} = \Im \left( \phi \overline{D_{\mu} \phi} \right)
  = \Im \left( \phi \overline{\partial_{\mu} \phi} \right)
  - A_{\mu} \abs{\phi}^{2},
\end{equation}
satisfying the conservation law
\begin{equation}\label{KGcurrentConserved}
  \partial^{\mu} j_{\mu} = 0.
\end{equation}
In fact, one has the general identity
$\partial_{\mu} \Im \left( \phi \overline{D^{\mu} \phi} \right)
= \Im \left( \phi \overline{D_{\mu} D^{\mu} \phi} \right)$,
so \eqref{KGcurrentConserved} follows immediately from \eqref{KG}.

The Maxwell-Klein-Gordon system is then obtained by coupling 
\eqref{KG} to the Maxwell equation
\begin{equation}\label{Maxwell}
  \partial^{\nu} F_{\mu \nu} = j_{\mu}
\end{equation}
where $F_{\mu\nu} = \partial_{\mu} A_{\nu} - \partial_{\nu} A_{\mu}$ 
is the electromagnetic field tensor and $j_{\mu}$ is the Klein-Gordon 
current density \eqref{KGcurrent}.
The system \eqref{Maxwell},\eqref{KGcurrent},\eqref{KG} is then 
what we --- provisionally --- call the Maxwell-Klein-Gordon system.
We want to consider this as a system of second order PDE in the 
unknowns $A_{\mu}$ and $\phi$, but there is an obvious problem with 
this, since $F_{\mu\nu}$ --- and hence the observables, i.e., the
electric and magnetic field 
vectors, whose components are entries of the matrix $F_{\mu\nu}$ ---
are not uniquely determined by $A_{\mu}$. This is known as the gauge 
ambiguity, and to resolve it one adds another equation to the system, 
a so-called gauge condition, which uniquely determines $A_{\mu}$.
The standard gauge conditions are
(i) Lorentz: $\partial^{\mu} A_{\mu} = 0$,
(ii) Coulomb: $\partial^{i} A_{i} = 0$ and
(iii) temporal: $A_{0} = 0$.

In this paper, we shall rely on the Coulomb condition, which carries 
the advantage --- as Klainerman and Machedon observed in 
\cite{Kl-Ma0.2} for the case of $n = 3$ --- that the bilinear terms involving
derivatives turn out to be of null form type, and therefore have better regularity 
properties than generic products.
Since the derivation of the null form structure in \cite{Kl-Ma0.2} uses the
special vector calculus of $n = 3$, in particular the $\curl$ operator,
we include a generalization of this argument to arbitrary dimension
in section \ref{Reformulation}.
\subsection{Main result}\label{Results}
If we add the Coulomb gauge condition $\partial^{j} A_{j} = 0$ to the MKG
system \eqref{Maxwell}, \eqref{KGcurrent}, \eqref{KG}
and expand, we get:
\begin{subequations}\label{MKG}
\begin{align}
   \label{MKGA}
   \Delta A_{0} &= - \Im \bigl(\phi \overline{\partial_t \phi} \bigr) +
   \abs{\phi}^2 A_0, \\
   \label{MKGB}
   \square A_{j} &= - \Im \bigl(\phi \overline{\partial_j \phi} \bigr) +
   \abs{\phi}^2 A_j - \partial_j \partial_t A_0,
   \\
   \label{MKGC}
   \square \phi &= - 2i A^j \partial_j \phi + 2i A_0 \partial_t \phi + i
   (\partial_t A_0) \phi + A^\mu A_\mu \phi + m^{2} \phi, \\
   \label{MKGD}
   \partial^j A_j &= 0.
\end{align}
\end{subequations}
In the rest of the paper, with the exception of section \ref{Scaling},
we will take $n = 4$. Thus, the unknowns are
$$
   A_0, A_j : \R^{1+4} \to \R,
   \quad
   \phi : \R^{1+4} \to \C.
$$
When convenient, we shall write $A$ for the four-vector field 
$(A^{j})_{j = 1,\dots,4}$. 
Initial data are specified at time $t = 0$:
\begin{subequations}\label{MKGData}
\begin{alignat}{2}
   \label{MKGDataA}
   A \init &= a \in H^s, \qquad & \partial_t A \init &= b \in H^{s-1}, \\
   \label{MKGDataB}
   \phi \init &= \phi_0 \in H^s, \qquad & \partial_t \phi \init &= 
\phi_1 \in H^{s-1},
\end{alignat}
\end{subequations}
where $H^s = \bigl\{ f \in \Schwartz'(\R^4) : (I - \Delta)^{s/2} f 
\in L^2(\R^4) \bigr\}$
and $a, b$ are real vector fields.
In view of the Coulomb condition \eqref{MKGD}, we must require
\begin{equation}\label{DivergenceFreeData}
   \partial^j a_j = \partial^j b_j = 0.
\end{equation}
Observe that no data are specified for the non-dynamical variable $A_{0}$. 
This is quite natural, because $A_{0}$ is determined by $\phi$ and 
$\partial_{t} \phi$ at any time $t$ by solving the elliptic equation 
\eqref{MKGA}.
\begin{theorem}\label{MKGTheorem} For all $s > 1$, the
Cauchy problem \eqref{MKG},\eqref{MKGData},\eqref{DivergenceFreeData}
on $\R^{1+4}$ is locally well-posed.
\end{theorem}
Local well-posedness here includes (a) existence of a local solution
\begin{subequations}\label{Reg}
\begin{align}
  \label{A0Reg}
   A_0 &\in C\bigl( [0,T], \dot H^1 \bigr) \cap C^{1}\bigl( [0,T], L^{2} \bigr)
   \\
   \label{AjPhiReg}
   A_j, \phi
   &\in C\bigl( [0,T], H^{s} \bigr) \cap C^{1}\bigl( [0,T], 
   H^{s-1} \bigr)
\end{align}
\end{subequations}
up to a time $T > 0$ depending continuously
on the $H^s$-norm of the initial data; (b) uniqueness of the solution; (c) 
continuous dependence on the data; and (d) persistence of higher regularity.
A more precise statement, for an equivalent system, can be found in 
Theorem \ref{WaveSystemTheorem}, section \ref{Reformulation}.
In particular, the uniqueness is proved not in the class \eqref{Reg}, 
but in a smaller space determined by the iteration norms; see 
\eqref{MKGRegularity}.

To prove Theorem \ref{MKGTheorem} we shall in effect eliminate the nondynamical
variable $A_0$ from the equations, by solving the elliptic equations.
This leaves us with a system of nonlinear wave equations,
which we then prove is locally well-posed. Once this has been achieved,
we can go back to the original system \eqref{MKG}, and conclude that this
is also well-posed.

Let us be more precise. We introduce a new variable
$B_0 = \partial_t A_0$. Applying $\partial^j$ to \eqref{MKGB}
and using \eqref{MKGD} yields
\begin{equation}\label{B0}
   \Delta B_0 = - \Im \partial^j \bigl(\phi \overline{\partial_j \phi} \bigr) +
   \partial^j \bigl( \abs{\phi}^2 A_j \bigr).
\end{equation}
Now we eliminate $A_0$ and $\partial_t A_0 = B_0$ from \eqref{MKGB} and
\eqref{MKGC} by solving \eqref{MKGA} and \eqref{B0}. Thus
$A_0 = A_0(\phi)$ and $B_0 = B_0(A,\phi)$ are nonlinear operators. Since
the Coulomb condition \eqref{MKGD} turns out to be automatically satisfied
because of the constraint
\eqref{DivergenceFreeData}, we obtain a system of nonlinear 
wave equations
\begin{subequations}\label{WaveSystem}
\begin{align}
   \label{WaveSystemA}
   \square A &= \mathcal M(A,\phi),
   \\
   \label{WaveSystemB}
   \square \phi &= \mathcal N(A,\phi),
\end{align}
\end{subequations}
where $\mathcal M$ and $\mathcal N$ are certain 
operators\footnote{See section \ref{Reformulation} for precise 
definitions}, nonlocal in the
space variable, which are sums of terms of the following types:
(i) bilinear and higher order multilinear expressions involving $A$ and $\phi$
and their first derivatives, (ii) terms involving $A_0(\phi)$, and
(iii) a linear term $m^{2} \phi$ in \eqref{WaveSystemB}. Moreover, all
the bilinear terms have a null structure, due to the Coulomb gauge,
and for these terms one already has good estimates (see
\cite{Kl-Ma5}, and also \cite{Kl-Ta} for the case of Yang-Mills; here we shall
rely more particularly on variants of these estimates proved in \cite{Kl-Se}).
We complement these with estimates for the higher
order multilinear terms and terms containing $A_0(\phi)$,
and local well-posedness of the system \eqref{WaveSystem} then
follows by the general theory developed in the author's paper \cite{Se}.

Then the original system \eqref{MKG} is also locally well-posed,
by reversing the steps leading to \eqref{WaveSystem}. That is, if $(A,\phi)$
has the requisite regularity (see \eqref{MKGRegularity}) and solves \eqref{WaveSystem}
on a time-slab, and if we set $A_0 = A_0(\phi)$, then $\partial_t A_0 = B_0(A,\phi)$
in the sense of distributions
and the triple $(A_0,A,\phi)$ solves \eqref{MKG} on the same time-slab.

Thus, we show that the systems \eqref{MKG} and
\eqref{WaveSystem} are \emph{equivalent} for sufficiently regular 
solutions.
\subsection{Scaling, optimality and the null condition}\label{Scaling}
As for many other field theories, there are two types of ``critical'' 
behaviour associated to the MKG system on $\R^{1+n}$. On the one hand, there is 
the critical regularity $s_{c}$ such that the homogeneous
initial data space $\dot H^{s_{c}}$ is left
invariant under the natural scaling transformation associated to MKG:
\begin{equation}\label{MKGscaling}
  A_{\mu}(t,x), \phi(t,x) \longrightarrow
  \lambda A_{\mu}(\lambda t,\lambda x), \lambda \phi(\lambda 
  t,\lambda x),
\end{equation}
where $\lambda$ is a positive parameter.\footnote{By this we mean that 
if $A_{\mu}, \phi$ solve MKG, then so do the rescaled fields, 
although the rest mass changes from $m$ to $\lambda m$.}
Since
\begin{equation}\label{DataScaling}
  \Sobdotnorm{\lambda f(\lambda x)}{s} = \lambda^{s-(n-2)/2} \Sobdotnorm{f}{s},
\end{equation}
we conclude that $s_{c} = \frac{n-2}{2}$. In general\footnote{See
\cite[Section 1.3]{Kl-Se} for further discussion and references.}
one expects field theories to be locally well-posed (LWP) for $s > 
s_{c}$ and ill-posed for $s < s_{c}$; we say more about this below.
In the critical case $s = s_{c}$ 
one expects some type of weakened well-posedness\footnote{For example, one does
not expect smooth dependence on initial data, which rules out proof by
iteration. A good example is wave maps into a sphere; see
Tao \cite[Section 1]{TaoWM1} for a summary of the regularity results 
for wave maps.} for data with small norm.

On the other hand, there is the energy-critical dimension $n$
such that the critical regularity is at the level of the 
energy:\footnote{MKG has a conserved energy which is at the level 
of the $H^{1}$ data norm; see \cite{Kl-Ma0.2}.}
$s_{c} = 1$. For MKG this means $n = 4$, which is the 
dimension we consider in this paper. For field theories in general, 
one expects global regularity in the critical dimension, as well as 
in subcritical dimensions ($s_{c} < 1$), and breakdown of regularity for
large data in supercritical dimensions ($s_{c} > 1$).

As mentioned above, the global regularity is known in the subcritical 
dimension $n = 3$ for MKG, but the question of global regularity in the 
critical dimension $n = 4$, even for data with small energy, remains open.
By conservation of energy, a LWP result, for small-norm data, at
the critical regularity $s_{c} = 1$ would settle this question in the 
affirmative, but it is perhaps more realistic to expect a more direct proof 
of global regularity in analogy with the results of Tao \cite{TaoWM1, TaoWM2}
for wave maps into a sphere. It is to be hoped that our almost optimal 
LWP result will play some role in any such result.

The expectation of ill-posedness for $s < s_{c}$ is based
on the scaling \eqref{MKGscaling} and 
\eqref{DataScaling}. First, if blow-up occurs for smooth, compactly 
supported data, then one can construct data in $H^{s}$, $s < s_{c}$, 
with arbitrarily small norm, for which 
there is no local existence; see, e.g., \cite[pp 98--99]{So} for
this argument. However, this is not a very convincing point to make 
here, as we do expect global regularity for MKG on $\R^{1+4}$.
We can show, however, that it is impossible to prove any 
well-posedness result for $s < s_{c}$ using an iteration argument based 
on estimates. The idea can be illustrated by the following example:
As is well-known, the algebra inequality
\begin{equation}\label{SobAlgebra}
  \Sobnorm{fg}{s} \le C_{s,n} \Sobnorm{f}{s} \Sobnorm{g}{s}
\end{equation}
holds for $H^{s}(\R^{n})$ iff $s > n/2$. A rather crude way of ruling out 
the range $s < n/2$ is to observe that if \eqref{SobAlgebra} holds, 
then by rescaling\footnote{In the limit $\lambda \to 
\infty$, the inhomogeneous Sobolev norm $H^{s}$ scales like $\dot H^{s}$.}
$x \to \lambda x$ and letting $\lambda \to \infty$, we get $1 \lesssim 
\lambda^{s-n/2}$.

This idea is easily applied to the iteration for MKG written in the 
form \eqref{WaveSystem}. Let us take $m = 0$ 
here to make the system scale invariant. If $a = b = 0$ in 
\eqref{MKGDataA} and $\phi_{1} = 0$ in \eqref{MKGDataB}, then
the first iterate of $A$ solves\footnote{Here $\Proj$ denotes the 
projection onto divergence free vector fields. See section 
\ref{Reformulation}.} $\square A^{(1)} = - 
\Proj \Im ( \phi^{(0)} \overline{ \nabla \phi^{(0)}})$ with zero data, 
where $\phi^{(0)}$ is the solution of $\square \phi^{(0)} = 0$ with 
data $(\phi_{0},0)$. If we can prove LWP in $H^{s}$ by 
iteration, there must be an estimate
\begin{equation}\label{FirstIterate}
  \sup_{0 \le t \le 1} \Sobnorm{A^{(1)}(t)}{s} \lesssim \Sobnorm{\phi_{0}}{s}^{2},
\end{equation}
for all $\phi_{0}$ with sufficiently small norm.
Now assume $s < s_{c}$. We then claim that \eqref{FirstIterate} implies
$A^{(1)} \equiv 0$, which is absurd.
Indeed, given $T > 0$, apply \eqref{FirstIterate} to the rescaled iterate
$$
  \widetilde A^{(1)} (t,x) = \lambda A^{(1)}(\lambda t ,\lambda x)
$$
at time $t = T/\lambda$. As $\lambda \to \infty$ this gives
$$
  \lambda^{s-s_{c}} \Sobnorm{A^{(1)}(T)}{s} \lesssim 
  (\lambda^{s-s_{c}})^{2} \Sobnorm{\phi_{0}}{s}^{2},
$$
whence $A^{(1)}(T) = 0$.
\begin{plainremark}
This argument has nothing to do with the null condition, of course. 
A more careful analysis (see \cite[Section 
1]{Kl-Se}) suggests that for a generic equation of the form
$$
  \square u = u \partial u
$$
on $\R^{1+n}$ one needs $s \ge \max(\tfrac{n-2}{2}, \tfrac{n+1}{4})$ in 
order for the iterates to stay in $H^{s}$, and this is consistent 
with Lindblad's counterexamples \cite{Li1}. However, if the right hand 
side is replaced by a null form expression like
\eqref{ProjectionFormula1} or \eqref{ProjectionFormula2}, then one 
only needs $s \ge \max(\tfrac{n-2}{2}, \tfrac{n-1}{4})$, so
the null condition improves matters when $n \le 4$.
\end{plainremark}
\medskip
As remarked already, the main difficulty is to prove LWP when $s$ is 
very close to $s_{c}$, whereas simpler arguments can be used for
larger $s$. Let us be more precise.
Observe that relative to Lorentz gauge, MKG on $\R^{1+n}$
is a system of nonlinear wave equations of the schematic form
(see \cite[Section 1]{Kl-Se})
\begin{equation}\label{MKGLorentz}
  \square u = u \partial u + u^{3},
\end{equation}
and for this system LWP for $s > n/2$ can be proved by standard methods,
just using the energy inequality for the wave equation and Sobolev embeddings.
This can easily be improved to $s > \tfrac{n-1}{2}$ by using a
$L_{t}^{2} L_{x}^{\infty}$ spacetime estimate instead of just Sobolev embedding.
For $n = 4$ this gives LWP for $s > 3/2$, which is still one quarter of a derivative
above what one expects (cf. remark above) from the analysis of the first iterate of
\eqref{MKGLorentz}, namely $s > 5/4$. No proof of LWP of \eqref{MKGLorentz} in this
range seems to exist in the literature, but it should be obtainable 
using the spaces $H^{s,\theta}$ (see section \ref{Notation}) and 
$L^{2}$ bilinear estimates for the homogeneous wave equation of the 
type first proved in \cite{Kl-Ma2}.
However, to go below the regularity $5/4$, one really needs the null condition, 
which seems to rule out Lorentz gauge. Of course, once a LWP result has been 
proved in one gauge, one can in principle use gauge transformations (see 
\cite{Kl-Ma0.2}) to transfer this result to other gauges; but to make this 
rigorous requires sufficient regularity of the solutions, and we will 
not consider this question here. 
\subsection{Function spaces}\label{Notation}
Here we define the spaces that we make use of. See \cite{Kl-Se} for 
more details.

The Fourier transform of $f(x)$ [resp. $u(t,x)$] is denoted
$\widehat f(\xi) = \Fourier f(\xi)$
[resp. $\widehat u(\tau,\xi) = \Fourier u(\tau,\xi)$].

We say that a norm $\norm{\cdot}$, on some space $\X$ of tempered 
distributions,
depends only on the size of the Fourier transform if
$$
   \abs{\widehat u} \le \abs{\widehat v} \implies \norm{u} \le \norm{v}.
$$
(Here we assume, of course, that the Fourier transform
of any element of $\X$ is a function.)

If $\X$ and $\Y$ are two normed spaces, the notation $\X \hookrightarrow \Y$
means continuous inclusion.

For any $\alpha \in \R$ we define Fourier multiplier operators
$\Lambda^{\alpha}$, $\Lambda_{+}^{\alpha}$ and $\Lambda_{-}^{\alpha}$ by
\begin{align*}
   \widehat{\Lambda^{\alpha} f}(\xi) &= \bigl( 1 + \abs{\xi}^{2}
   \bigr)^{\alpha/2} \widehat{f}(\xi),
   \\
   \widehat{\Lambda_{+}^{\alpha} u}(\tau,\xi) &= \bigl( 1 + \tau^{2} + 
\abs{\xi}^2
   \bigr)^{\alpha/2} \widehat{u}(\tau,\xi),
   \\
   \widehat{\Lambda_{-}^{\alpha} u}(\tau,\xi) &=
   \left( 1 + \frac{(\tau^2 - \abs{\xi}^2)^{2} }
   {1 + \tau^{2} + \abs{\xi}^2 } \right)^{\alpha/2}
    \widehat{u}(\tau,\xi).
\end{align*}
It should be remarked that the weight of $\Lambda_{-}^{\alpha}$
is comparable to $\bigl(1 + \hypwt{\tau}{\xi}\bigr)^{\alpha}$,
but the former has the advantage of being smooth.

The Sobolev and ``Wave Sobolev'' spaces $H^s$ and $H^{s,\theta}$ are given
by the weighted $L^{2}$ norms
$$
  \Sobnorm{f}{s} = \twonorm{\Lambda^s f}{(\R^4)}
  \quad \text{and}
  \quad
  \spacetimenorm{u}{s}{\theta} = \twonorm{\Lambda^s \Lambda_-^\theta 
  u}{(\R^{1+4})}.
$$
We shall also use the related space $\scrH^{s,\theta}$ defined by
$$
   \Spacetimenorm{u}{s}{\theta} = \spacetimenorm{u}{s}{\theta}
   + \spacetimenorm{\partial_t u}{s-1}{\theta}
   \sim \twonorm{\Lambda^{s-1} \Lambda_+ \Lambda_-^\theta u}{}.
$$
In view of Plancherel's theorem, these norms depend only on the size 
of the Fourier transform.
It is an important fact that when $\theta > 1/2$, the spaces $H^{s,
\theta}$ and $\scrH^{s,\theta}$ can be localized in time, since then the
embeddings
\begin{equation}\label{EnergyEmbedding}
   H^{s,\theta} \hookrightarrow C_b(\R,H^s) \quad \text{and} \quad
   \scrH^{s,\theta} \hookrightarrow C_b(\R,H^s) \cap C^1_b(\R,H^{s-1})
\end{equation}
hold. See \cite[Section 3]{Kl-Se}.

Since $L^2( \abs{\xi}^2 d\xi ) \subseteq L^1_{\text{loc}}(\R^4)
\subseteq \Schwartz'(\R^4)$, we may define
$$
   \dot H^1 = \Fourier^{-1} \bigl[ L^2( \abs{\xi}^2 d\xi ) \bigr].
$$
Thus $\dot H^1$ is a Hilbert space
with norm $\Sobdotnorm{f}{1}^2 = \int_{\R^4} \abs{\xi}^2
\bigabs{\widehat f(\xi)}^2 \, d\xi$.
We remark that if $\dot W^1 = \{ f : \nabla f \in L^2 \}$, then
$\dot H^1$ is obtained by identifying elements of $\dot W^1$
differing by a constant. Observe also that $\Schwartz$ is dense in
$L^2( \abs{\xi}^2 d\xi )$, hence in $\dot H^1$.
We shall use frequently the fact that
\begin{equation}\label{L4Sob}
   \dot H^1 \hookrightarrow L^4(\R^4).
\end{equation}
In other words, $\Lpnorm{f}{4} \lesssim \Sobdotnorm{f}{1}$.
This holds by the Hardy-Littlewood-Sobolev inequality (see Stein 
\cite[Chapter V]{St}).

If $\X$ is a separable Banach space of functions on $\R^4$, and $1 
\le p \le \infty$,
we denote by $L_t^p(\X)$ the space $L^p(\R,\X)$ of $\X$-valued functions.
In particular, we write
$$
   \mixednorm{u}{p}{q} = \left( \int_\R \norm{u(t,\cdot)}_{L^q(\R^4)}^p
   \, dt \right)^{1/p}
$$
with the usual modification if $p = \infty$.

We also need a version of this last norm which only depends on the size
of the Fourier transform: If $u \in \Schwartz'$ and $\widehat u$ is a tempered
function, set
$$
   \Mixednorm{u}{p}{q} = \sup \left\{ \int_{\R^{1+4}} \abs{ \widehat u (\tau,\xi)
   } \widehat v(\tau,\xi) \, d\tau \, d\xi : \text{$v \in \Schwartz$, $\widehat v \ge 0$,
   $\mixednorm{v}{p'}{q'} = 1$} \right\},
$$
where $1 = \frac{1}{p} + \frac{1}{p'}$ and $1 = \frac{1}{q} + \frac{1}{q'}$.
Let $\Mixed{p}{q}$ be the corresponding subspace of $\Schwartz'$.
Then $\Mixednorm{\cdot}{p}{q}$ is a translation invariant norm on
$\Mixed{p}{q}$. Note that $\Mixed{2}{2} = L^2(\R^{1+4})$ and
\begin{equation}\label{ModifiedMixednormBound}
   \Mixednorm{u}{p}{q} \le \mixednorm{u}{p}{q} \quad \text{whenever} \quad
   \widehat u \ge 0.
\end{equation}
We refer the reader to \cite[Section 4]{Kl-Se} for more details on 
these spaces.

We can now make precise the regularity statement \eqref{Reg}.
The solutions we obtain are in the following spaces:
\begin{subequations}\label{MKGRegularity}
\begin{align}
   A_0 &\in C\bigl( [0,T], \dot H^1 \bigr) \cap C^{1}\bigl( [0,T], L^{2} \bigr),
   \\
   \label{MKGRegularityAj}
   A_j &\in \scrH^{s,\theta} \cap
   \Lambda^{-\gamma} \Lambda_-^{-\half} \bigl[ \Mixed{1}{8} \bigr],
   \\
   \label{MKGRegularityphi}
   \phi &\in \scrH^{s,\theta},
\end{align}
\end{subequations}
where $\theta > \half$ and $\gamma > 0$ depend on $s$. For technical 
reasons, it is useful to iterate $A_{j}$ and $\phi$ in these global 
spaces, but in the end we are only interested in their values on a
time interval $[0,T]$ whose size depends on the norms of the data. Since 
the space $\scrH^{s,\theta}$ can be localized in time, this presents 
no problems.

\begin{remark}\label{SpacesRemark}
  The auxiliary space $\Mixed{1}{8}$ in \eqref{MKGRegularityAj}
  is necessary when $s < 5/4$. See Theorem 8.2 in \cite{Kl-Se} and the remark 
  following it.
\end{remark}

\noindent
\emph{Note:}
Throughout the paper, we use the convenient shorthand $\lesssim$
for $\le$ up to a positive multiplicative constant $C$.
Usually $C$ is completely innocuous, and only depends
on parameters that may be considered fixed. There are exceptions,
notably for Lipschitz estimates (then $C$ is only ``locally'' constant),
but these are clearly pointed out.
\subsection{Reformulation of the MKG system}\label{Reformulation}
As discussed in section \ref{Results}, an important step in our proof is
to recast the MKG system \eqref{MKG} as a system of nonlinear
wave equations \eqref{WaveSystem}. Here we describe this in detail.

As was shown in \cite{Kl-Ma0.2}, the first terms
on the right hand sides of equations \eqref{MKGB} and \eqref{MKGC}
can be expressed, due to the Coulomb condition \eqref{MKGD},
in terms of the bilinear null forms
\begin{equation}\label{Qij}
   Q_{jk}(u,v) = \partial_j u \, \partial_k v - \partial_k u \, \partial_j v.
\end{equation}
Since the argument in \cite{Kl-Ma0.2} was special to the case $n = 3$,
we include here a proof of this fact which works for any dimension.
First, let $\Proj$ be the projection onto the divergence free 
vector fields on $\R^4$. In terms of the Riesz transforms
$R_j = (-\Delta)^{-\half} \partial_j$,
$$
   \Proj X_j = X_j + R_j R^k X_k = R^k ( R_j X_k  - R_k X_j).
$$
Observe that $\Proj$ is bounded on every $L^p$, $1 < p < \infty$,
since this is true for the Riesz transforms (see Stein \cite{St}).
Moreover, it is clear that the Riesz transforms, and hence
$\Proj$, are bounded on any space whose norm only depends on the size 
of the Fourier transform, in particular on any Sobolev space $H^s$.

Since $\partial_{j} (u \partial_{k} v) - \partial_{k} (u \partial_{j} v)
= Q_{jk}(u,v)$, it follows immediately from the definition of $\Proj$ 
that
\begin{equation}\label{ReformulationA}
   \Proj (u \partial_j v) = R^k (-\Delta)^{-\half} Q_{jk} (u,v),
\end{equation}
whence
\begin{equation}\label{ProjectionFormula1}
   \Proj \left( - \Im \bigl[ \phi \overline{\partial_j \phi} \bigr] \right)
   = 2 R^k (-\Delta)^{-\half} Q_{jk} (\Re \phi, \Im \phi).
\end{equation}
Also,
$$
   2 \partial_j u \Proj X^j = Q_{jk} \left( u,
   (-\Delta)^{-\half} \bigl[ R^j X^k - R^k X^j \bigr] \right),
$$
as one can see by expanding the right hand side.
Therefore, if $A$ is divergence free, so that $\Proj A = A$, then
\begin{equation}\label{ProjectionFormula2}
   2 A^j \partial_j \phi
   = Q_{jk} \left( \phi,
   (-\Delta)^{-\half} \bigl[ R^j A^k - R^k A^j \bigr] \right).
\end{equation}
\begin{remark}\label{RegRemark}
The calculations leading to the identity \eqref{ReformulationA} are certainly
justified when $u$ and $v$ belong to the Schwartz class $\Schwartz(\R^4)$.
Moreover, both sides of the identity are bounded bilinear
operators of $(u,v) \in H^s \times H^s$ into $H^{-1}$, where $s > 1$. 
Thus the identity holds
for all $u,v \in H^s$, and we conclude that \eqref{ProjectionFormula1}
holds for all $\phi$ with the regularity \eqref{MKGRegularityphi},
since by \eqref{EnergyEmbedding} this implies $\phi \in C_b(\R,H^s)$.
To bound the left hand side of \eqref{ReformulationA}, use first the 
dual
\begin{equation}\label{DualL4Sob}
   \bigtwonorm{(-\Delta)^{-\half} f}{(\R^4)} \lesssim 
   \norm{f}_{L^{4/3}(\R^4)}
\end{equation}
of \eqref{L4Sob}. Since $\Proj$ is bounded on $L^{p}$,
it then suffices to observe that
\begin{equation}\label{RegRemark1}
   \Lpnorm{u \partial_j v}{4/3} \lesssim
   \Lpnorm{u}{4} \Lpnorm{\partial_j u}{2}
   \lesssim \Sobnorm{u}{1} \Sobnorm{v}{1},
\end{equation}
where we used \eqref{L4Sob}.
To prove boundedness of the right hand side of \eqref{ReformulationA},
it is enough to show
$$
   \bigSobnorm{(-\Delta)^{-\half} (fg)}{-1}
   \lesssim \Sobnorm{f}{s-1}\Sobnorm{g}{s-1}.
$$
This can be reduced, via the self-duality of $L^2$,
Plancherel's theorem, and the Cauchy-Schwarz inequality, to the fact that
$\abs{\xi}^{-1} (1 + \abs{\xi})^{-1-2(s-1)}$ belongs to $L^2(\R^4)$,
since $s > 1$.
Similar, but simpler, considerations show that the remaining
bilinear and cubic terms in (\ref{MKG}a,b,c) and (\ref{MKG'}a,c,d)
are bounded into $C_b\bigl(\R,L^{4/3}(\R^4)\bigr)$ when regarded as
operators on $A_0,A,\phi$ in the class \eqref{MKGRegularity}.
For example, for a cubic expression $uvw$
we have by H\"older's inequality and \eqref{L4Sob}
that
\begin{equation}\label{RegRemark2}
  \Lpnorm{uvw}{4/3}
  \le \Lpnorm{u}{4} \Lpnorm{v}{4} \Lpnorm{w}{4}
  \lesssim
  \Sobdotnorm{u}{1}\Sobdotnorm{v}{1}\Sobdotnorm{w}{1}.
\end{equation}
\end{remark}
\medskip
Returning to the main thread of our argument, we now use the null 
form identities derived above
to arrive at an equivalent formulation of MKG:
\begin{subequations}\label{MKG'}
\begin{align}
   \label{MKG'A}
   \Delta A_{0} &= - \Im \bigl(\phi \overline{\partial_t \phi} \bigr) +
   \abs{\phi}^2 A_0, \\
   \label{MKG'B}
   \Delta \partial_t A_{0}
   &= - \Im \partial^j \bigl(\phi \overline{\partial_j \phi} \bigr) +
   \partial^j \bigl( \abs{\phi}^2 A_j \bigr)
   \\
   \label{MKG'C}
   \square A_{j} &= 2 R^k (-\Delta)^{-\half} Q_{jk} (\Re \phi, \Im \phi)
   + \Proj \bigl( \abs{\phi}^2 A_j \bigr)
   \\
   \label{MKG'D}
   \square \phi &= - i Q_{jk} \left( \phi,
   (-\Delta)^{-\half} \bigl[ R^j A^k - R^k A^j \bigr] \right)
   \\
   \notag
   &\quad + 2i A_0 \partial_t \phi + i
   (\partial_t A_0) \phi + A^\mu A_\mu \phi + m^{2} \phi.
\end{align}
\end{subequations}
This system acts as a stepping stone between \eqref{MKG} and 
\eqref{WaveSystem}.
\begin{proposition}\label{ReformulationPropositionA}
The systems \eqref{MKG} and \eqref{MKG'} are equivalent.
More precisely, any local solution of \eqref{MKG} with the
regularity \eqref{MKGRegularity} and divergence free initial
data is a solution of \eqref{MKG'} and vice versa.
\end{proposition}
\begin{proof}
To go from \eqref{MKG} to \eqref{MKG'}, observe that $A_j$ is divergence
free by \eqref{MKGD}; apply $\partial^j$ to \eqref{MKGB} to
get equation \eqref{MKG'B}; apply $\Proj$ to \eqref{MKGB} and
use \eqref{ProjectionFormula1} to get \eqref{MKG'C}; finally,
\eqref{MKG'D} follows from \eqref{MKGC} using \eqref{ProjectionFormula2}.

To go the other way, observe that by \eqref{ProjectionFormula1},
the right hand side of \eqref{MKG'C} is divergence free;
thus $\square \partial^j A_j = 0$, and since the initial data
of $A_j$ are divergence free, \eqref{MKGD} follows. Then, in
view of \eqref{ProjectionFormula2}, \eqref{MKGC}
and \eqref{MKG'D} are equivalent. Finally, to go from \eqref{MKG'C}
to \eqref{MKGB}, it suffices to check that the right hand side
of the latter is divergence free. But this follows from \eqref{MKG'B}.
\end{proof}
Once the system has been written in the form \eqref{MKG'}
it is easy to eliminate $A_0$ and $\partial_t A_0$ and
obtain the system of wave equations \eqref{WaveSystem}.
We now describe this in more detail.
\begin{lemma}\label{A0Lemma}
Given $\phi$ in the class \eqref{AjPhiReg},
equation \eqref{MKG'A} has a unique solution $A_0 \in \dot H^1$
on every time-slice $\{t\} \times \R^4$, and these solutions
assemble to a space-time function $A_0 = A_0(\phi) \in C_b(\R,\dot H^1)$.
Moreover, we have bounds, on every time-slice $\{t\} \times \R^4$,
$$
   \Sobdotnorm{A_0}{1} \le 2 \twonorm{\partial_t \phi}{}
$$
and
$$
   \Sobdotnorm{A_0(\phi)-A_0(\psi)}{1} \lesssim
   \Sobnorm{\phi-\psi}{1}
   + \twonorm{\partial_t \phi - \partial_t \psi}{},
$$
where the suppressed constant depends polynomially on
$\Sobnorm{\phi}{1}, \Sobnorm{\psi}{1},
\twonorm{\partial_t \phi}{}$ and $\twonorm{\partial_t \psi}{}$,
but is independent of $t$.
\end{lemma}
This is proved in section \ref{EllipticEstimates}.

Next we consider \eqref{MKG'B}, with $\partial_t A_0$
replaced by the new variable $B_0$:
\begin{equation}\label{MKG'BB}
   \Delta B_{0}
   = - \Im \partial^j \bigl(\phi \overline{\partial_j \phi} \bigr) +
   \partial^j \bigl( \abs{\phi}^2 A_j \bigr).
\end{equation}
\begin{lemma}\label{B0Lemma}
Given $(A,\phi)$ in the class \emph{(\ref{AjPhiReg})},
the equation \eqref{MKG'BB} has a unique solution $B_0 \in L^2$
on every time-slice $\{t\} \times \R^4$, given by
\begin{equation}\label{B0Def}
   B_0 = R^j (- \Delta)^{-\half} \left[
   \Im \bigl(\phi \overline{\partial_j \phi} \bigr) -
   \abs{\phi}^2 A_j \right],
\end{equation}
and the solutions assemble to a space-time function
$B_0 = B_0(A,\phi) \in C_b(\R,L^2)$. Moreover,
we have bounds, on every time-slice $\{t\} \times \R^4$,
$$
   \twonorm{B_0}{}
   \le C \bigl( 1 + \Sobnorm{A}{1} \bigr) \Sobnorm{\phi}{1}^2
$$
for a constant $C$ independent of $t$, and
$$
   \twonorm{B_0(A,\phi)-B_0(A',\phi')}{}
   \lesssim \Sobnorm{A-A'}{1} + \Sobnorm{\phi - \phi'}{1},
$$
where the suppressed constant depends polynomially on
$\Sobnorm{A}{1}, \Sobnorm{A'}{1},
\Sobnorm{\phi}{1}$ and $\Sobnorm{\phi'}{1}$, but is independent of $t$.
\end{lemma}
\begin{proof} To see that \eqref{B0Def} is in $L^2$,
first apply \eqref{DualL4Sob}, then estimate as in \eqref{RegRemark1} 
and \eqref{RegRemark2}. That \eqref{B0Def} is the only $L^2$ solution
can be seen by taking the Fourier transform of both sides of \eqref{MKG'BB}.
\end{proof}
In view of the above lemmas, \eqref{MKG'} implies \eqref{WaveSystem}, with
$$
   \mathcal M = (\mathcal M_1,\dots,\mathcal M_4),
   \quad \mathcal M_j = \mathcal M_{j,1} + \mathcal M_{j,2},
   \quad
   \mathcal N = \mathcal N_1 + \cdots + \mathcal N_6,
$$
where
\begin{align*}
   \mathcal M_{j,1} &=
   2 R^k (-\Delta)^{-\half} Q_{jk} (\Re \phi, \Im \phi),
   \\
   \mathcal M_{j,2} &=
   \Proj \bigl( \abs{\phi}^2 A_j \bigr),
   \\
   \mathcal N_1 &= - i Q_{jk} \left( \phi,
   (-\Delta)^{-\half} \bigl[ R^j A^k - R^k A^j \bigr] \right),
   \\
   \mathcal N_2 &=
   2i A_0(\phi) \partial_t \phi,
   \\
   \mathcal N_3 &=
   i B_0(A,\phi) \phi,
   \\
   \mathcal N_4 &=
   - [A_0(\phi)]^2 \phi,
   \\
   \mathcal N_5 &=
   \abs{A}^2 \phi,
   \\
   \mathcal N_6 &= m^{2} \phi,
\end{align*}
and $\abs{A}^2 = A_j A^j$ in the next to last line.

Arguing as in Remark \ref{RegRemark},
and using Lemmas \ref{A0Lemma} and \ref{B0Lemma}, it is readily checked
that the multilinear expressions in $\mathcal M$ and $\mathcal N$ are 
all continuous
maps into $C_b(\R,L^{4/3})$ [or $C_b(\R,H^{-1})$ in the case of 
$\mathcal M_{j,1}$]
for $(A,\phi)$ in the class (\ref{MKGRegularity}b,c). However, proving
the following theorem
requires much more sophisticated estimates.
\begin{theorem}\label{WaveSystemTheorem}
The system of wave equations \eqref{WaveSystem}, with $\mathcal M$
and $\mathcal N$ defined as above, is locally well-posed for initial
data in $H^s$, all $s > 1$, in the following sense (all pairs $(A,\phi)$
are understood to belong to the class \emph{(\ref{MKGRegularity}b,c)}
in what follows):
\begin{itemize}
  \item[(a)] {\bfseries(Local existence)}
   For all initial data \eqref{MKGData} there exists a $T > 0$,
   which depends continuously on the norms of the data, and
   there exists a pair $(A,\phi)$ which solves \eqref{WaveSystem} in
   the sense of distributions
   on $(0,T) \times \R^4$ and satisfies the given initial condition.
   \item[(b)] {\bfseries(Uniqueness)} If $T > 0$ and we have two
   solutions $(A,\phi)$ and $(A',\phi')$ of \eqref{WaveSystem} on
   $(0,T) \times \R^4$ with identical initial data, then they agree on 
the entire
   time-slab.
   \item[(c)] {\bfseries(Continuous dependence on initial data)}
   If, for some $T > 0$, $(A,\phi)$ solves \eqref{WaveSystem} on
   $(0,T) \times \R^4$ with initial data \eqref{MKGData}, then for
   all initial data $(a',b',\phi_0',\phi_1')$ such that
   $$
     \delta = \Sobnorm{a-a'}{s} + \Sobnorm{b-b'}{s-1}
     + \Sobnorm{\phi_0 - \phi_0'}{s} + \Sobnorm{\phi_1 - \phi_1'}{s-1}
   $$
   is sufficiently small, there is a solution $(A',\phi')$ on the
   same time-slab and with these initial data. Moreover, we have
   $$
     \Sobnorm{A-A'}{s} + \Sobnorm{\partial_t A-\partial_tA'}{s-1}
     + \Sobnorm{\phi-\phi'}{s} + \Sobnorm{\partial_t \phi-\partial_t \phi'}{s-1}
     \le C \delta
   $$
   uniformly in $0 \le t \le T$.
   \item[(d)]
   {\bfseries(Persistence of higher regularity)}
   If $k$ is a positive integer and $(A,\phi)$ solves
   \eqref{WaveSystem} on $(0,T) \times \R^4$ with initial data
   in $H^{s+k}$ (that is, \eqref{MKGData} holds with $s$ replaced by
   $s+k$), then
   $$
     A,\phi \in C\bigl([0,T],H^{s+k}\bigr) \cap
     C^{1}\bigl([0,T],H^{s+k-1}\bigr).
   $$
   \item[(e)] {\bfseries (Classical solutions)}
   If the data belong to $H^{s+k}$ for every $k$, then
   the solution is smooth:
   $$
     A,\phi \in C^\infty\bigl( [0,T] \times \R^4 \bigr).
   $$
\end{itemize}
\end{theorem}
The proof of this theorem will occupy us in the next two sections.

Here we want to show that Theorem \ref{MKGTheorem}
can be deduced from Theorem \ref{WaveSystemTheorem}.
It clearly suffices to demonstrate the equivalence of the
systems \eqref{MKG} and \eqref{WaveSystem}.
The remainder of this section is devoted to a proof of this fact,
assuming that the conclusions of Theorem \ref{WaveSystemTheorem} hold.
\begin{proposition}\label{ReformulationPropositionB}
The systems \eqref{MKG} and \eqref{WaveSystem}
are equivalent for local solutions in the regularity
class \eqref{MKGRegularity}, with divergence free initial data.
\end{proposition}
In view of Proposition \ref{ReformulationPropositionA}, it
suffices to show the equivalence of \eqref{MKG'} and \eqref{WaveSystem}.
We have seen already that \eqref{MKG'} implies \eqref{WaveSystem}.
The converse is not quite so obvious, but for sufficiently regular solutions
it follows by some straightforward calulations and
the fact, proved in section \ref{EllipticEstimates}, that the
only $\dot H^1$ solution of the elliptic equation
$\Delta u = \abs{\phi}^2 u$ is $u = 0$.
For general $H^s$ data we then choose an approximating
sequence of sufficiently regular data, use the persistence of higher regularity
and continuous dependence on initial data, which hold by virtue of
Theorem \ref{WaveSystemTheorem}, and pass to the limit.

We now turn to the details.

Assume that $(A,\phi)$ is in the
class (\ref{MKGRegularity}b,c) and solves \eqref{WaveSystem} on a time-slab
$S_T = (0,T) \times \R^4$, with initial data satisfying
\eqref{MKGData} and \eqref{DivergenceFreeData}. Set $A_0 = A_0(\phi)$.
Then \eqref{MKG'} is satisfied, but with $\partial_t A_0$ replaced
by $B_0 = B_0(A,\phi)$
in \eqref{MKG'B} and \eqref{MKG'D}.
Thus, all we have to prove is that the distributional derivative
$\partial_t A_0$ agrees with $B_0$ on $S_T$. At first glance one may 
think that this is simply a matter of taking a time derivative of 
\eqref{MKG'A} and using the conservation law \eqref{KGcurrentConserved} 
to conclude that $\Delta \partial_{t} A_{0} = \Delta B_{0}$, but this 
is a circular argument since the derivation of 
\eqref{KGcurrentConserved} is not valid unless we know that 
$\partial_{t} A_{0} = B_{0}$.

In what follows, keep in mind that $A_\mu$ and $B_0$ are real-valued.\
Applying $\partial_t$ to \eqref{MKG'A} gives
\begin{equation}\label{ReformulationD}
   \Delta \partial_t A_0 = - \Im \bigl(\phi \overline{\partial_t^2 
\phi} \bigr) +
   2 \Re \bigl( \phi \overline{\partial_t \phi} \bigr) A_0
   + \abs{\phi}^2 \partial_t A_0.
\end{equation}
Since \eqref{MKG'C} and \eqref{DivergenceFreeData} hold,
it follows as in the proof of Proposition
\ref{ReformulationPropositionA} that $A$ is divergence free.
Therefore, \eqref{ProjectionFormula2} holds, and
since \eqref{MKG'D} holds (with $\partial_t A_0$ replaced by $B_0$),
we conclude that
$$
   -\partial_t^2 \phi + \Delta \phi
   = \square \phi = - 2i A^j \partial_j \phi + 2i A_0 \partial_t \phi
   + i B_0 \phi + A^\mu A_\mu \phi + m^{2} \phi.
$$
Using this expression for $\partial_t^2 \phi$ gives, after some calculation,
$$
   - \Im \bigl(\phi \overline{\partial_t^2 \phi} \bigr)
   = - \Im \partial^j \bigl(\phi \overline{\partial_j \phi} \bigr)
   + 2 \Re \bigl( \phi \overline{\partial_j \phi} \bigr) A^j
   - 2 \Re \bigl( \phi \overline{\partial_t \phi} \bigr) A_0
   - B_0 \abs{\phi}^2.
$$
Since
$$
   - \Im \partial^j \bigl(\phi \overline{\partial_j \phi} \bigr)
   = \Delta B_0 - \partial^j \bigl( \abs{\phi}^2 A_j \bigr),
$$
we get
$$
   - \Im \bigl(\phi \overline{\partial_t^2 \phi} \bigr)
   = \Delta B_0 - \partial^j \bigl( \abs{\phi}^2 A_j \bigr)
   + 2 \Re \bigl( \phi \overline{\partial_j \phi} \bigr) A^j
   - 2 \Re \bigl( \phi \overline{\partial_t \phi} \bigr) A_0
   - B_0 \abs{\phi}^2.
$$
Inserting this in \eqref{ReformulationD} gives
$$
   \Delta \partial_t A_0 =
   \Delta B_0 - \partial^j \bigl( \abs{\phi}^2 A_j \bigr)
   + 2 \Re \bigl( \phi \overline{\partial_j \phi} \bigr) A^j
   - B_0 \abs{\phi}^2
   + \abs{\phi}^2 \partial_t A_0.
$$
But
$$
   \partial^j \bigl( \abs{\phi}^2 A_j \bigr)
   = 2 \Re \bigl( \phi \overline{\partial_j \phi} \bigr) A^j
   + \abs{\phi}^2 \partial^j A_j
   = 2 \Re \bigl( \phi \overline{\partial_j \phi} \bigr) A^j
$$
since $A$ is divergence free, and so we finally get
$$
   \Delta (\partial_t A_0 - B_0) =
   \abs{\phi}^2 (\partial_t A_0 - B_0).
$$
The above manipulations are justified provided
\begin{equation}\label{HigherReg1}
   \partial_t A_0 \in C([0,T],\dot H^1).
\end{equation}
If, moreover,
\begin{equation}\label{HigherReg2}
   B_0 \in C([0,T],\dot H^1),
\end{equation}
then it follows by the uniqueness result alluded to above
(see Lemma \ref{EllipticLemma} in section \ref{EllipticEstimates})
that $\partial_t A_0 = B_0$ in $[0,T] \times \R^4$.

But \eqref{HigherReg1} and \eqref{HigherReg2} certainly hold
under the additional assumption that the initial data
\eqref{MKGData} of $A$ and $\phi$ belong to
$H^{s+k}$ for every positive integer $k$.
Leaving aside the proof of this assertion for the moment,
we note that any $f \in H^s$ can be approximated in the $H^s$ norm
by a sequence belonging to every $H^{s+k}$, by convolution with a
$C_c^\infty$ approximation of the identity, and if $f$ is
divergence free, then so is the approximating sequence.
Combining these facts with the continuous dependence of $A$ and $\phi$
on their $H^s$ initial data (Theorem \ref{WaveSystemTheorem}),
and the continuity of the operators $A_0$ and $B_0$ (Lemmas \ref{A0Lemma}
and \ref{B0Lemma}), we conclude by passing to the limit that the equality
$\partial_t A_0 = B_0$ holds in the sense of distributions
on $(0,T) \times \R^4$ for all initial data \eqref{MKGData}
satisfying \eqref{DivergenceFreeData}.

It remains to prove that \eqref{HigherReg1} and \eqref{HigherReg2}
hold if the initial data \eqref{MKGData} of $A$ and $\phi$ belong to
$H^{s+k}$ for every positive integer $k$.
For $A_0$, this follows by persistence of higher regularity
(part (d) of Theorem \ref{WaveSystemTheorem}),
the inductive regularity step \eqref{ClassicalA} in
section \ref{ClassicalSolutions} and
Lemma \ref{SpaceTimeEllipticLemma} in the same section. As for $B_0$, in view of
\eqref{B0Def} it is clear that, on every time-slice,
$$
  \Sobdotnorm{B_0}{1} \le \sum_j \left( \twonorm{\phi \partial_j \phi}{}
  + \bigtwonorm{\abs{\phi}^2 A_j}{} \right)
$$
and by H\"older's inequality and Sobolev embedding it is
easy to see that the right hand side is dominated by
$\Sobnorm{\phi}{1} \Sobnorm{\phi}{2} + \Sobnorm{\phi}{2}^2 \Sobnorm{A}{2}$.
But if $A$ and $\phi$ have initial data in $H^{s+1}$, then
by persistence of higher regularity
(part (d) of Theorem \ref{WaveSystemTheorem})
we know that $A,\phi \in C([0,T],H^2)$.
\section{Proof of Theorem \ref{WaveSystemTheorem}}\label{MainEstimates}
Here we discuss the estimates needed to prove local
well-posedness of the system \eqref{WaveSystem}, with
$\mathcal M$ and $\mathcal N$ defined as in section
\ref{Reformulation}.

The local existence for the system \eqref{WaveSystem} is proved 
by Picard iteration in the spaces \eqref{MKGRegularityAj} and 
\eqref{MKGRegularityphi}, 
which are defined using the spacetime Fourier transform, and hence are 
global. However, since they embed in \eqref{AjPhiReg}, they can easily
be localized in time. In fact, this time localization smooths out the
singularity of the inverse $\square^{-1}$ of the wave operator, and
--- if done with sufficient care --- allows one to handle large initial
data by taking a sufficiently small time interval. These matters are
considered in detail in the author's paper \cite{Se}, and also in \cite[Section 
5]{Kl-Se}, and we refer the interested reader there.

Fix $1 < s < 2$. (For larger $s$, the result can be proved
by simpler arguments.) Let $\theta > \half$ and $\gamma, \varepsilon > 0$;
these quantities depend on the choice of $s$, and will be specified
later. Now define
\begin{align*}
   \X_1 &= \scrH^{s,\theta} \cap
   \Lambda^{-\gamma} \Lambda_-^{-\half} \bigl[ \Mixed{1}{8} \bigr],
   \\
   \X_2 &= \scrH^{s,\theta},
   \\
   \Y_k &= \Lambda_+ \Lambda_-^{1-\varepsilon} \X_k, \quad k = 1,2
\end{align*}
with norms
\begin{align*}
   \norm{A}_{\X_1} &= \Spacetimenorm{A}{s}{\theta}
   + \bigMixednorm{\Lambda^\gamma \Lambda_-^{\half} A}{1}{8},
   \\
   \norm{\phi}_{\X_2} &= \Spacetimenorm{\phi}{s}{\theta},
   \\
   \norm{F}_{\Y_k} &= \norm{\Lambda_+^{-1} \Lambda_-^{-1+\varepsilon} 
F }_{\X_k},
   \quad k = 1,2.
\end{align*}
All these spaces are complete (see \cite[Proposition 4.2]{Kl-Se}),
and by \cite[Proposition 5.6]{Kl-Se}, $\X_1$ and $\X_2$ satisfy the
hypotheses of \cite[Theorem 1]{Se}. Consequently, by \cite[Theorem 2]{Se},
the system \eqref{WaveSystem} is locally well-posed for $H^s$ data
if the following Lipschitz conditions\footnote{Keep in mind
that $\mathcal M$ and $\mathcal N$ vanish
at the origin, so if we take $A' = 0$ and $\phi' = 0$, we
simply get bounds for $\mathcal M(A,\phi)$ and $\mathcal N(A,\phi)$.}
hold:
\begin{subequations}\label{LipschitzEstimates}
\begin{align}
   \norm{\mathcal M(A,\phi) - \mathcal M(A',\phi')}_{\Y_1}
   &\lesssim \norm{A-A'}_{\X_1} + \norm{\phi-\phi'}_{\X_2},
   \\
   \norm{\mathcal N(A,\phi) - \mathcal N(A',\phi')}_{\Y_2}
   &\lesssim \norm{A-A'}_{\X_1} + \norm{\phi-\phi'}_{\X_2},
\end{align}
\end{subequations}
where the suppressed constants depend continuously on
$$
   \norm{A}_{\X_1}, \quad \norm{A'}_{\X_1},
   \quad \norm{\phi}_{\X_2} \quad \text{and} \quad \norm{\phi'}_{\X_2}.
$$
In fact, these estimates
guarantee that the conclusions (a,b,c) of Theorem \ref{WaveSystemTheorem}
hold. In the next section we show how to prove parts (d) and (e) of 
the same theorem.

It suffices to prove \eqref{LipschitzEstimates} with $\mathcal M$
replaced by $\mathcal M_{j,k}$ and with $\mathcal N$ replaced by
$\mathcal N_1,\dots,\mathcal N_5$. Furthermore, in view of the multilinear
structure, it suffices to prove (concerning the 
suppressed constants, see note below):
\begin{align}
   \label{M1}
   \norm{\mathcal M_{j,1}}_{\Y_1} &\lesssim \norm{\phi}_{\X_2}^2,
   \\
   \label{M2}
   \norm{\mathcal M_{j,2}}_{\Y_1} &\lesssim \norm{A}_{\X_1} 
\norm{\phi}_{\X_2}^2,
   \\
   \label{N1}
   \norm{\mathcal N_1}_{\Y_2} &\lesssim \norm{A}_{\X_1} \norm{\phi}_{\X_2},
   \\
   \label{N2}
   \norm{\mathcal N_2}_{\Y_2} &\lesssim \norm{A_0(\phi)}_{\Z_1} 
\norm{\phi}_{\X_2},
   \\
   \label{N3}
   \norm{\mathcal N_3}_{\Y_2} &\lesssim \norm{B_0(A,\phi)}_{\Z_2} 
\norm{\phi}_{\X_2},
   \\
   \label{N4}
   \norm{\mathcal N_4}_{\Y_2} &\lesssim \norm{A_0(\phi)}_{L_t^\infty(\dot H^1)}
   \norm{A_0(\phi)}_{\Z_1} \norm{\phi}_{\X_2},
   \\
   \label{N5}
   \norm{\mathcal N_5}_{\Y_2} &\lesssim \norm{A}_{\X_1}^2 \norm{\phi}_{\X_2},
   \\
   \label{N6}
   \norm{\mathcal N_6}_{\Y_2} &\le \norm{\phi}_{\X_2},
   \\
   \label{A0EstA}
   \norm{A_0}_{L_t^\infty(\dot H^1)} &\lesssim \norm{\phi}_{\X_2},
   \\
   \label{A0EstAA}
   \norm{A_0(\phi)-A_0(\phi')}_{L_t^\infty(\dot H^1)}
   &\lesssim \norm{\phi-\phi'}_{\X_2},
\end{align}
where $\Z_1$ and $\Z_2$ are certain intermediate spaces, to be
specified later, such that
\begin{align}
   \label{A0EstB}
   \norm{A_0}_{\Z_1} &\lesssim
   \norm{\phi}_{\X_2}^2 + \norm{\phi}_{\X_2}^3,
   \\
   \label{A0EstBB}
   \norm{A_0(\phi)-A_0(\phi')}_{\Z_1} &\lesssim \norm{\phi-\phi'}_{\X_2},
   \\
   \label{B0EstA}
   \norm{B_0}_{\Z_2} &\lesssim \bigl( 1 + \norm{A}_{\X_1} \bigr) 
\norm{\phi}_{\X_2}^2,
   \\
   \label{B0EstB}
   \norm{B_0(A,\phi) - B_0(A',\phi')}_{\Z_2} &\lesssim
   \norm{A-A'}_{\X_1} + \norm{\phi-\phi'}_{\X_2}.
\end{align}
It should be emphasized that in the Lipschitz estimates \eqref{A0EstAA}, \eqref{A0EstBB}
and \eqref{B0EstB}, the suppressed constant depends polynomially on
the norms $\norm{\phi}_{\X_2}$ and $\norm{\phi'}_{\X_2}$,
and in the case of \eqref{B0EstB} also on $\norm{A}_{\X_1}$ and 
$\norm{A'}_{\X_1}$. Observe that the estimate \eqref{N6} for the 
linear term is trivial, since the norms only depend on the size of the 
Fourier transform.

The following was proved in \cite[Theorem 8.6]{Kl-Se}.
\begin{plaintheorem}
The estimates \eqref{M1} and \eqref{N1} hold provided
\begin{subequations}\label{Parameters}
\begin{gather}
   \label{Parameters1}
   \half < \theta < \min \left( \frac{3}{4}, \frac{s}{2}\right)
   \\
   \label{Parameters2}
   0 < \varepsilon < \frac{1}{4} \min
   \left( \frac{3}{4} - \theta, \frac{s}{2} - \theta \right)
   \\
   \gamma = \theta - \half - 3\varepsilon.
\end{gather}
\end{subequations}
\end{plaintheorem}
Having fixed $\theta$ and $\varepsilon$ satisfying these
requirements, we define $p$ and $r$ by
\begin{equation}\label{prdef}
   \frac{1}{p} = \frac{3}{2} - \theta - 2\varepsilon,
   \quad \frac{1}{r} = 1 - \theta - 2\varepsilon,
\end{equation}
and we choose $q$ so large that
\begin{equation}\label{qdef}
   \frac{4}{q} < \min \left( 2\theta - 1, 1-\frac{1}{p} \right).
\end{equation}
Observe that as $s \to 1$, the triple $(p,q,r) \to (1,\infty,2)$. Now set
\begin{align}
   \label{Z1}
   \norm{A_0}_{\Z_1} &= \mixednorm{\Lambda^{s-1} A_0}{p}{q},
   \\
   \label{Z2}
   \norm{B_0}_{\Z_2} &= \mixednorm{\Lambda^{s-1} B_0}{r}{8/3}.
\end{align}

For easy reference, we list here some estimates that we shall use
(here $p,q,r$ are defined as above):
\begin{align}
   \label{Tataru}
   \bigmixednorm{\Lambda^{s-1}(-\Delta)^{-1}(uv)}{p}{q}
   &\lesssim \spacetimenorm{u}{s}{\theta} \spacetimenorm{v}{s-1}{\theta},
   \\
   \label{Pecher0}
   \mixednorm{u}{2}{8} &\lesssim \spacetimenorm{u}{1}{\theta},
   \\
   \label{Pecher1}
   \mixednorm{u}{r}{8} &\lesssim \spacetimenorm{u}{s}{\theta},
   \\
   \label{Pecher2}
   \mixednorm{u}{2p}{\beta} &\lesssim \spacetimenorm{u}{s}{\theta},
   \quad \frac{5}{2} + \theta + 2\varepsilon - 2s \le \frac{8}{\beta} 
\le 2\theta,
   \\
   \label{L1L2Embedding}
   \spacetimenorm{u}{0}{\theta+\varepsilon-1} &\lesssim \mixednorm{u}{p}{2},
   \\
   \label{DualEmbedding}
   \Mixednorm{u}{1}{8} &\lesssim
   \bigMixednorm{\Lambda \Lambda_-^{\half+\varepsilon} u}{1}{2},
   \\
   \label{CalculusIneq}
   \Sobnorm{fg}{\sigma} &\lesssim \Lpnorm{\Lambda^\sigma f}{p_1} \Lpnorm{g}{q_1}
   + \Lpnorm{\Lambda^\sigma g}{p_2} \Lpnorm{f}{q_2},
\end{align}
where in the last inequality,
$$
   \sigma > 0, \quad \frac{1}{p_k} + \frac{1}{q_k} = \half,
   \quad 2 \le p_k < \infty.
$$

The inequality \eqref{Tataru} follows from a theorem of 
Klainerman-Tataru \cite{Kl-Ta};
we give the details in an appendix.

The Strichartz type estimates (\ref{Pecher0}--\ref{Pecher2}) are 
special cases of
\cite[Theorem D]{Kl-Se}. (The [non-optimal] upper bound for $8/\beta$ 
in \eqref{Pecher2}
guarantees that the pair $(2p,\beta)$ is wave admissible; the lower bound
is chosen so that we do not exceed $s$ space derivatives on the right 
hand side.)

The inequality \eqref{L1L2Embedding} can either be proved directly, using
Plancherel's theorem, H\"older's inequality, Minkowski's integral inequality
and the Hausdorff-Young inequality, or it can be proved by interpolation,
as in \cite[Section 6(vii)]{Kl-Se}).

Inequality \eqref{DualEmbedding}
is a special case of \cite[Proposition 4.8]{Kl-Se}.

The calculus inequality \eqref{CalculusIneq} is Lemma 1 in Ponce-Sideris
\cite{Po-Si}.

As mentioned already, \eqref{M1} and \eqref{N1} hold by
\cite[Theorem 8.6]{Kl-Se}. We now prove the remaining estimates
\eqref{M2} and (\ref{N2}--\ref{B0EstB}), thereby concluding
the proof of parts (a,b,c) of Theorem \ref{WaveSystemTheorem}.
\subsection{Proof of \eqref{M2}}
Since the norm only depends on the size of the Fourier transform,
we can ignore the projection $\mathcal P$. More accurately,
$$
   \norm{\mathcal M_{j,2}}_{\Y_1} \lesssim \bignorm{ \abs{\phi}^2 A_j}_{\Y_1}.
$$
Thus, it suffices to prove
$$
   \norm{\Lambda_+^{-1} \Lambda_-^{\varepsilon-1} (uvw) }_{\X_1}
   \lesssim \Spacetimenorm{u}{s}{\theta}
   \Spacetimenorm{v}{s}{\theta} \Spacetimenorm{w}{s}{\theta},
$$
or, equivalently,
\begin{align*}
   \spacetimenorm{uvw}{s-1}{\theta+\varepsilon-1}
   &\lesssim \Spacetimenorm{u}{s}{\theta}
   \Spacetimenorm{v}{s}{\theta} \Spacetimenorm{w}{s}{\theta},
   \\
   \bigMixednorm{\Lambda^\gamma \Lambda_+^{-1} 
\Lambda_-^{\varepsilon-\half} (uvw) }{1}{8}
   &\lesssim \Spacetimenorm{u}{s}{\theta}
   \Spacetimenorm{v}{s}{\theta} \Spacetimenorm{w}{s}{\theta}.
\end{align*}
Since all the norms depend only on the size of the Fourier transform,
we may assume that $u,v,w$ have non-negative Fourier transforms, and
we see that it is sufficient to prove
(note that $\gamma + 2\varepsilon < s-1$ by \eqref{Parameters})
\begin{align}
   \label{M2A}
   \spacetimenorm{uvw}{0}{\theta+\varepsilon-1}
   &\lesssim \Spacetimenorm{u}{1}{\theta}
   \Spacetimenorm{v}{s}{\theta} \Spacetimenorm{w}{s}{\theta},
   \\
   \label{M2B}
   \bigMixednorm{\Lambda^{-1} \Lambda_-^{-\varepsilon-\half} (uvw) }{1}{8}
   &\lesssim \Spacetimenorm{u}{1}{\theta}
   \Spacetimenorm{v}{s}{\theta} \Spacetimenorm{w}{s}{\theta}.
\end{align}
By \eqref{L1L2Embedding} and H\"older's inequality,
$$
   \spacetimenorm{uvw}{0}{\theta+\varepsilon-1}
   \lesssim \mixednorm{u}{\infty}{4} \mixednorm{v}{2p}{8} \mixednorm{w}{2p}{8},
$$
and \eqref{M2A} follows by Sobolev embedding and \eqref{Pecher2}.

Using \eqref{DualEmbedding} and \eqref{ModifiedMixednormBound}, we get
$$
   \bigMixednorm{\Lambda^{-1} \Lambda_-^{-\varepsilon-\half} (uvw) }{1}{8}
   \lesssim \mixednorm{uvw}{1}{2}
   \lesssim \mixednorm{u}{\infty}{4} \mixednorm{v}{2}{8} \mixednorm{w}{2}{8}.
$$
Now use Sobolev embedding and \eqref{Pecher0}.
\subsection{Proof of \eqref{N2}}
We have to show
$$
   \spacetimenorm{uv}{s-1}{\theta+\varepsilon-1}
   \lesssim \mixednorm{\Lambda^{s-1} u}{p}{q} \Spacetimenorm{v}{s-1}{\theta}.
$$
By \eqref{L1L2Embedding} and \eqref{CalculusIneq},
\begin{multline*}
   \spacetimenorm{uv}{s-1}{\theta+\varepsilon-1}
   \lesssim \mixednorm{\Lambda^{s-1}(uv)}{p}{2}
   \\
   \lesssim \mixednorm{\Lambda^{s-1} u}{p}{q}
\mixednorm{v}{\infty}{(1/2-1/q)^{-1}}
   + \mixednorm{u}{p}{\infty} \mixednorm{\Lambda^{s-1} v}{\infty}{2}.
\end{multline*}
The desired estimate now follows by Sobolev embedding, since 
$\frac{4}{q} < s-1$.
\subsection{Proof of \eqref{N3}}
We must prove
$$
   \spacetimenorm{uv}{s-1}{\theta+\varepsilon-1}
   \lesssim \mixednorm{\Lambda^{s-1} u}{r}{8/3} \Spacetimenorm{v}{s}{\theta}.
$$
By \eqref{L1L2Embedding} and \eqref{CalculusIneq},
\begin{multline}\label{N3ProofA}
   \spacetimenorm{uv}{s-1}{\theta+\varepsilon-1}
   \lesssim \mixednorm{\Lambda^{s-1}(uv)}{p}{2}
   \\
   \lesssim \mixednorm{\Lambda^{s-1} u}{2}{8/3} \mixednorm{v}{r}{8}
   + \mixednorm{u}{r}{8/3} \mixednorm{\Lambda^{s-1} v}{2}{8}.
\end{multline}
Now apply \eqref{Pecher0} and \eqref{Pecher1}. Note also that
$\mixednorm{u}{r}{8/3} \lesssim \mixednorm{\Lambda^{s-1} u}{r}{8/3}$,
since $\Lambda^{-\delta}$ is bounded on $L^p$ for all $1 \le p \le \infty$
and $\delta \ge 0$. In fact, $\Lambda^{-\delta}$ corresponds
to convolution with an $L^1$ function; see Stein \cite{St}.
\subsection{Proof of \eqref{N4}}
It suffices to show
$$
   \spacetimenorm{u^2 v}{s-1}{\theta+\varepsilon-1}
   \lesssim \mixednorm{u}{\infty}{4}
   \mixednorm{\Lambda^{s-1} u}{p}{q} \Spacetimenorm{v}{s}{\theta}.
$$
By \eqref{L1L2Embedding} and \eqref{CalculusIneq},
\begin{align*}
   &\spacetimenorm{u^2 v}{s-1}{\theta+\varepsilon-1}
   \lesssim \mixednorm{\Lambda^{s-1}(u^2 v)}{p}{2}
   \\
   &\qquad \lesssim \mixednorm{\Lambda^{s-1} u}{p}{q} 
\mixednorm{uv}{\infty}{(1/2-1/q)^{-1}}
   + \mixednorm{u}{\infty}{4} \mixednorm{\Lambda^{s-1} (uv)}{p}{4}
   \\
   &\qquad \lesssim \mixednorm{\Lambda^{s-1} u}{p}{q} \mixednorm{u}{\infty}{4}
   \mixednorm{v}{\infty}{(1/4-1/q)^{-1}}
   \\
   &\qquad \qquad \qquad \qquad \qquad \qquad \qquad \quad
   + \mixednorm{u}{\infty}{4} \mixednorm{u}{p}{\infty}
   \mixednorm{\Lambda^{s-1} v}{\infty}{4}.
\end{align*}
Now apply Sobolev embedding, and use \eqref{qdef}.
\subsection{Proof of \eqref{N5}}
It suffices to show
$$
   \spacetimenorm{u v w}{s-1}{\theta+\varepsilon-1}
   \lesssim \Spacetimenorm{u}{s}{\theta}
   \Spacetimenorm{v}{s}{\theta} \Spacetimenorm{w}{s}{\theta},
$$
but this was proved above; see the proof of \eqref{M2}.
\subsection{Proof of \eqref{A0EstA} and \eqref{A0EstAA}}
These follow from Lemma \ref{A0Lemma}, which is proved
in section \ref{EllipticEstimates}.
\subsection{Proof of \eqref{A0EstB} and \eqref{A0EstBB}}
Since
$$
   A_0 = (-\Delta)^{-1} \left[ \Im \bigl(\phi \overline{\partial_t \phi} \bigr)
   - \abs{\phi}^2 A_0 \right],
$$
it suffices, taking into account the multilinearity of the terms
inside the brackets, as well as the estimates \eqref{A0EstA} and 
\eqref{A0EstAA},
to show that
\begin{align*}
   \mixednorm{\Lambda^{s-1} (-\Delta)^{-1} (uv)}{p}{q} &\lesssim
   \spacetimenorm{u}{s}{\theta} \spacetimenorm{v}{s-1}{\theta},
   \\
   \mixednorm{\Lambda^{s-1} (-\Delta)^{-1} (uvw)}{p}{q} &\lesssim
   \spacetimenorm{u}{s}{\theta} \spacetimenorm{v}{s}{\theta}
   \mixednorm{w}{\infty}{4}.
\end{align*}
The former is exactly \eqref{Tataru}, and the left hand side of the 
latter is $\lesssim$
\begin{equation}\label{A0EstB1}
   \mixednorm{(-\Delta)^{-1} (uvw)}{p}{q}
   + \bigmixednorm{(-\Delta)^{\frac{s-3}{2}} (uvw)}{p}{q}.
\end{equation}
Here we applied the following useful result, which is an immediate consequence
of Lemma 2(ii) in Chapter V of Stein \cite{St}.
\begin{lemma}\label{BesselRieszLemma}
For $\alpha > 0$ and $1 \le p \le \infty$,
$$
   \Lpnorm{\Lambda^\alpha f}{p} \lesssim \Lpnorm{f}{p} + 
\bigLpnorm{(-\Delta)^{\alpha/2} f}{p},
$$
where the suppressed constant only depends on $\alpha$.
\end{lemma}
Returning to the sum \eqref{A0EstB1}, note that by Sobolev embedding, 
it is $\lesssim$
$$
   \mixednorm{uvw}{p}{\alpha_1}
   + \bigmixednorm{uvw}{p}{\alpha_2},
$$
where
\begin{align*}
   \frac{1}{\alpha_1} &= \half + \frac{1}{q} = \frac{1}{4}
   + 2 \left( \frac{1}{8} + \frac{1}{2q} \right),
   \\
   \frac{1}{\alpha_2} &= \frac{3-s}{4} + \frac{1}{q} = \frac{1}{4}
   + 2 \left( \frac{2-s}{8} + \frac{1}{2q} \right).
\end{align*}
Thus
$$
   \mixednorm{uvw}{p}{\alpha_k} \le \mixednorm{u}{2p}{\beta_k}
   \mixednorm{v}{2p}{\beta_k} \mixednorm{w}{\infty}{4}, \quad k = 1,2
$$
where
$$
   \frac{1}{\beta_1} = \frac{1}{8} + \frac{1}{2q},
   \quad \frac{1}{\beta_2} = \frac{2-s}{8} + \frac{1}{2q}.
$$
Using \eqref{Parameters2} and \eqref{qdef} it is easily checked that
$$
   \frac{5}{2} + \theta + 2\varepsilon - 2s
   \le 2-s \le \frac{8}{\beta_2} < \frac{8}{\beta_1} < 2\theta,
$$
so we may apply \eqref{Pecher2} to finish the proof.
\subsection{Proof of \eqref{B0EstA} and \eqref{B0EstB}}\label{B0EstAProof}
We prove \eqref{B0EstA}; the same proof gives \eqref{B0EstB}
if one exploits the multilinearity of the terms
defining $B_0$.

First observe that by Lemma \ref{BesselRieszLemma},
$$
   \mixednorm{\Lambda^{s-1} B_0 }{r}{8/3}
   \lesssim \mixednorm{ B_0 }{r}{8/3}
   + \bigmixednorm{(-\Delta)^{\frac{s-1}{2}} B_0 }{r}{8/3}.
$$
Therefore, by Sobolev embedding, we have to estimate
$$
   \bigmixednorm{(-\Delta)^\half B_0 }{r}{\alpha_k}, \quad k = 1,2
$$
where
$$
   \frac{1}{\alpha_1} = \frac{5}{8}, \quad \frac{1}{\alpha_2} = 
\frac{5}{8} - \frac{s-1}{4}.
$$
Since $B_0$ is given by \eqref{B0Def}, and since the Riesz transforms
$R_j$ are bounded on $L^p$, $1 < p < \infty$, we see that it is enough to prove
\begin{align*}
   \mixednorm{uv}{r}{\alpha_k} &\lesssim
   \spacetimenorm{u}{s}{\theta} \spacetimenorm{v}{s-1}{\theta},
   \\
   \mixednorm{uvw}{r}{\alpha_k} &\lesssim
   \spacetimenorm{u}{s}{\theta} \spacetimenorm{v}{s}{\theta} 
\spacetimenorm{w}{s}{\theta}.
\end{align*}
By H\"older's inequality,
\begin{align*}
   \mixednorm{uv}{r}{\alpha_k} &\le
   \mixednorm{u}{r}{8} \mixednorm{v}{\infty}{\beta_k},
   \\
   \mixednorm{uvw}{r}{\alpha_k} &\le
   \mixednorm{u}{r}{8} \mixednorm{v}{\infty}{4} \mixednorm{v}{\infty}{\gamma_k},
\end{align*}
where
\begin{alignat*}{2}
   \beta_1 &= 2, &\qquad \gamma_1 &= 4,
   \\
   \frac{1}{\beta_2} &= \frac{1}{2} - \frac{s-1}{4},
   &\qquad \frac{1}{\gamma_2} &= \frac{1}{4} - \frac{s-1}{4}.
\end{alignat*}
Now apply \eqref{Pecher1} and Sobolev embedding.
\section{Higher regularity}
Here we prove parts (d) and (e) of Theorem \ref{WaveSystemTheorem}.
\subsection{The persistence property}
The key to proving part (d) of Theorem \ref{WaveSystemTheorem}
is to establish, for $k = 0,1,2,\dots$,
\begin{subequations}\label{HigherOrderEstimates}
\begin{align}
   \norm{\Lambda^k \mathcal M(A,\phi)}_{\Y_1}
   &\le \alpha_k
   \left\{ \norm{\Lambda^k A}_{\X_1} +  \norm{\Lambda^k \phi}_{\X_2} \right\}
   + \beta_k,
   \\
   \norm{\Lambda^k \mathcal N(A,\phi)}_{\Y_2}
   &\le \alpha_k
   \left\{ \norm{\Lambda^k A}_{\X_1} +  \norm{\Lambda^k \phi}_{\X_2} \right\}
   + \beta_k,
\end{align}
\end{subequations}
where
\begin{itemize}
   \item $\alpha_k$ depends continuously on
   $\norm{A}_{\X_1}$ and $\norm{\phi}_{\X_2}$,
   \item $\beta_0 = 0$,
   \item $\beta_k$, for $k \ge 1$, depends continuously on
   $\norm{\Lambda^{k-1} A}_{\X_1}$ and $\norm{\Lambda^{k-1} \phi}_{\X_2}$.
\end{itemize}
The case $k = 0$ is of course true by \eqref{LipschitzEstimates},
but it is useful to include it here for technical reasons.

In the absence of the lower order term $\beta_k$,
we could now appeal directly to \cite[Theorem 2]{Se},
to conclude that part (d) of Theorem \ref{WaveSystemTheorem}
holds. However, we can easily modify the proof given in \cite{Se}
to cover this more general case, as we demonstrate below.

First, however, let us dispose of proof of the above estimates.
Observe that we have the equivalence of norms
$$
   \norm{\Lambda^k u}_{\X_j}
   \sim \sum_{\abs{\alpha} \le k} \norm{\partial_x^\alpha u}_{\X_j}.
$$
This is trivial in view of the fact that the norms only depend
on the size of the Fourier transform. It is therefore clear, from the
multilinear structure of $\mathcal M$ and $\mathcal N$, and the product
rule for derivatives, that \eqref{HigherOrderEstimates} follows
from the very estimates proved in section \ref{MainEstimates}.
The only exception is the nonlinear operator $A_0(\phi)$,
for which we need the following estimate, replacing \eqref{A0EstA}:
\begin{lemma}\label{SpaceDerivativesEllipticLemma}
If $\Lambda^k \phi \in \X_2$, then
$$
   \norm{\partial_x^\alpha A_0}_{L_t^\infty(\dot H^1)}
   \le \gamma_k\left( \norm{\phi}_{\X_2} \right)
   \norm{\Lambda^k \phi}_{\X_2} + \eta_k\bigl( \norm{\Lambda^{k-1} 
  \phi}_{\X_2} \bigr) \quad \text{for all} \quad
   \abs{\alpha} \le k,
$$
where $\gamma_k$ and $\eta_k$ are continuous functions.
\end{lemma}
This is proved in section \ref{EllipticLemmaProofs}.

Let us now turn to the proof of Theorem \ref{WaveSystemTheorem}, part (d).

The issue is to show that if we have a pair $(A,\phi)$,
belonging to the class (\ref{MKGRegularity}b,c),
which solves \eqref{WaveSystem} on $S_T = (0,T) \times \R^4$
with initial data \eqref{MKGData}, and if the data have
some additional regularity, say $H^{s+k}$, then this extra
regularity persists throughout the time interval $[0,T]$:
\begin{equation}\label{HigherRegularitySpace}
   A,\phi \in C\bigl([0,T],H^{s+k}\bigr) \cap
   C^{1}\bigl([0,T],H^{s+k-1}\bigr).
\end{equation}
Now, as proved in \cite[Section 6.4]{Se}, it suffices to prove
this for \emph{some $T > 0$ which depends continuously on}
$$
   E_0 = \Sobnorm{a}{s} + \Sobnorm{b}{s-1}
   + \Sobnorm{\phi_0}{s} + \Sobnorm{\phi_1}{s-1}.
$$
We shall prove this using the Picard iterates corresponding to
the given data. It will be convenient to introduce the notation
$$
   E_k = \Sobnorm{a}{s+k} + \Sobnorm{b}{s+k-1}
   + \Sobnorm{\phi_0}{s+k} + \Sobnorm{\phi_1}{s+k-1}.
$$

Now fix an integer $K \ge 1$, and denote by $\alpha$ and $\beta$
the pointwise maxima of $\alpha_k$ and $\beta_k$, respectively,
taken over all $0 \le k \le K$. Let us assume that the initial
data belong to $H^{s+K}$, that is,
$$
   E_K < \infty.
$$

It is proved in \cite{Se} that for any $0 < T < 1$, there is a
linear operator $W_T$, which is bounded from $\Y_j \to \X_j$
($j = 1,2$), and such that $u = W_T F$ solves the inhomogeneous
wave equation $\square u = F$ on $(0,T) \times \R^4$
with vanishing initial data at $t = 0$.
Moreover, if $C_T$ is the maximum of the operator norms, that is,
\begin{equation}\label{OperatorNorm}
   C_T = \max \left( \norm{W_T}_{\Y_1 \to \X_1}, \norm{W_T}_{\Y_2 \to 
\X_2} \right),
\end{equation}
then
\begin{equation}\label{CTDecay}
   C_T \to 0 \quad \text{as} \quad T \to 0.
\end{equation}

The sequence of Picard iterates $(A^{(m)},\phi^{(m)})$ is
then defined inductively
as follows. First, let $A^{(0)}$ and $\phi^{(0)}$ be the solutions of
$\square A^{(0)} = 0$ and $\square \phi^{(0)} = 0$
with initial data \eqref{MKGData}, and then multiply them
by a smooth bump function which equals $1$ on the
interval $[0,T]$. By \cite[Theorem 1]{Se},
\begin{equation}\label{FirstIterateBound}
   \bignorm{\Lambda^k A^{(0)}}_{\X_1}
   + \bignorm{\Lambda^k \phi^{(0)}}_{\X_2} \le C E_k,
\end{equation}
with $E_k$ as above. Then define
\begin{align*}
   A^{(m+1)} &= A^{(0)} + W_T \mathcal M(A^{(m)},\phi^{(m)}),
   \\
   \phi^{(m+1)} &= \phi^{(0)} + W_T \mathcal N(A^{(m)},\phi^{(m)}).
\end{align*}
Let us write
\begin{align*}
   R^{(m)}_k &= \bignorm{\Lambda^k A^{(m)}}_{\X_1}
   + \bignorm{\Lambda^k \phi^{(m)}}_{\X_2},
   \\
   \omega^{(m)} &= \bignorm{A^{(m)} - A^{(m-1)}}_{\X_1}
   + \bignorm{\phi^{(m)} - \phi^{(m-1)}}_{\X_2}.
\end{align*}

Then by \eqref{FirstIterateBound}, \eqref{OperatorNorm} and
\eqref{HigherOrderEstimates} (with $k = 0$), we have
$$
   R^{(m+1)}_0 \le C E_0 + C_T \, \alpha\bigl( R^{(m)}_0 \bigr) R^{(m)}_0,
   \quad m \ge 0.
$$
If we choose $T$ so small that
\begin{equation}\label{Tcondition}
   2 C_T \, \alpha(2CE_0) \le 1,
\end{equation}
then it follows by induction on $m$ that
\begin{equation}\label{R0Bound}
   R^{(m)}_0 \le 2CE_0, \quad m \ge 0.
\end{equation}
Then, using the Lipschitz estimates \eqref{LipschitzEstimates}
(and making $\alpha$ larger if necessary),
$$
   \omega^{(m+1)} \le \half \omega^{(m)},
$$
so the sequence of Picard iterates is Cauchy in $\X_1 \times \X_2$,
and therefore converges; the limit
is of course the unique solution $(A,\phi)$ of our equation.

We shall prove that, \emph{with $T$ as in \eqref{Tcondition},}
\begin{equation}\label{Rkbound}
   R^{(m)}_k \le C_k(E_0,\dots,E_k), \quad k \le K, \quad m \ge 0,
\end{equation}
where $C_k$ is some continuous function.

Let us first see why this implies the desired conclusion
\eqref{HigherRegularitySpace} for $k \le K$.
The point is that by \eqref{Rkbound}, the sequence of
Picard iterates is bounded in the Hilbert space
$\scrH^{s+k,\theta}$ (recall that $\X_1 \hookrightarrow
\X_2 = \scrH^{s,\theta}$), and therefore, some
subsequence converges weakly in that space.
Since weak convergence in $\scrH^{s+k,\theta}$ implies
convergence in the sense of distributions, we conclude
that the strong limit $(A,\phi)$ agrees, as a distribution,
with this weak limit. Thus, $(A,\phi)$ belongs to
$\scrH^{s+k,\theta}$, and this immediately gives
\eqref{HigherRegularitySpace}.

We shall prove \eqref{Rkbound} by induction on $k$.

We already have the case $k = 0$, by \eqref{R0Bound}.

Now assume that $k < K$ and that \eqref{Rkbound}
holds. We claim that this implies \eqref{Rkbound} for $k+1$.
Indeed, by \eqref{HigherOrderEstimates},
\eqref{OperatorNorm} and \eqref{FirstIterateBound},
$$
   R^{(m+1)}_{k+1} \le C E_{k + 1}
   + C_T \alpha \bigl( R^{(m)}_0 \bigr) R^{(m)}_{k + 1}
   + C_T \beta\bigl( R^{(m)}_{k} \bigr).
$$
Taking into account \eqref{R0Bound}, \eqref{Tcondition}
and the induction hypothesis, we get
$$
    R^{(m+1)}_{k+1} \le C E_{k + 1}
   + \half R^{(m)}_{k + 1}
   + \frac{\beta\bigl( C_{k} (E_0,\dots,E_{k}) \bigr)}{2\alpha(2CE_0)}
$$
for $m \ge 0$. It now follows by induction on $m$ that
$$
    R^{(m)}_{k+1} \le 2C E_{k + 1}
   + \frac{\beta\bigl( C_{k} (E_0,\dots,E_{k}) \bigr)}{\alpha(2CE_0)},
   \quad m \ge 0,
$$
using \eqref{FirstIterateBound} for the case $m = 0$.
\subsection{Classical solutions}\label{ClassicalSolutions}
Here we outline the proof of part (e) of Theorem \ref{WaveSystemTheorem}.
In view of part (d) of the same theorem, it suffices to prove
the inductive step
\begin{equation}\label{ClassicalA}
   A,\phi \in \bigcap_{k = 1}^\infty C^m\bigl([0,T],H^{s+k}\bigr)
   \implies A,\phi \in \bigcap_{k = 1}^\infty C^{m+1}\bigl([0,T],H^{s+k}\bigr).
\end{equation}
But since $(A,\phi)$ solves \eqref{WaveSystem} on $(0,T) \times \R^4$,
we have there
\begin{align*}
   \partial_t^2 A &= \Delta A - \mathcal M(A,\phi),
   \\
   \partial_t^2 \phi &= \Delta \phi - \mathcal N(A,\phi),
\end{align*}
and so it is clear that \eqref{ClassicalA} follows from
\begin{multline}\label{ClassicalB}
   A,\phi \in \bigcap_{k = 1}^\infty C^m\bigl([0,T],H^{s+k}\bigr)
   \\
   \implies
   \mathcal M(A,\phi), \mathcal N(A,\phi)
   \in \bigcap_{k = 1}^\infty C^{m-1}\bigl([0,T],H^{s+k}\bigr).
\end{multline}

The key observation is of course that $\mathcal M$ and $\mathcal N$
only contain first order derivatives in time. Recall that $\mathcal M$
and $\mathcal N$ are sums of multilinear expressions in $A$ and $\phi$
and their first order derivatives, and terms involving $A_0(\phi)$.
But $A_0(\phi)$ is determined by the elliptic equation \eqref{MKG'A},
which also contains only first order partial derivatives in time
of $\phi$.

Thus, to prove \eqref{ClassicalB}, simply apply up to $m-1$ time
derivatives and any number of space derivatives, say $K$, to $\mathcal M$
and $\mathcal N$, and use the product rule for derivatives.
It is then easy to show --- we omit the details --- that on each time-slice,
the $L^2$-norms of the resulting expressions are bounded
in terms of (here $\alpha$ is a multi-index)
$$
  \bigSobnorm{\partial_t^j A}{K+k},
  \quad \bigSobnorm{\partial_t^j \phi}{K+k}
  \quad \text{and} \quad
  \bigSobdotnorm{\partial_t^j \partial_x^\alpha A_0(\phi)}{1}
$$
for $j \le m$, $\abs{\alpha} \le K$ and $k$ sufficiently large.
Then one appeals to the following higher regularity result for
$A_0(\phi)$, which is proved in section \ref{EllipticLemmaProofs}.
\begin{lemma}\label{SpaceTimeEllipticLemma}
Let $m, M$ be non-negative integers. If $\phi \in C^{m+1}\bigl([0,T],H^{M}\bigr)$,
that is, if
$$
   \partial_t^j \partial_x^\alpha \phi \in  C\bigl([0,T],L^2\bigr)
   \quad \text{for all} \quad j \le m+1 \quad \text{and all} \quad 
\abs{\alpha} \le M+1,
$$
where $\alpha$ is a multi-index, then
$$
   \partial_t^j \partial_x^\alpha A_0(\phi) \in  C\bigl([0,T],\dot H^1\bigr)
   \quad \text{for all} \quad j \le m \quad \text{and all} \quad 
\abs{\alpha} \le M,
$$
and $\bignorm{\partial_t^j \partial_x^\alpha A_0(\phi)}_{L^\infty([0,T],\dot H^1)}$
is bounded by a continuous function of the norms
$\norm{\partial_t^k\phi}_{L^\infty([0,T],H^{M})}$ for $k \le m+1$.
\end{lemma}
\section{Elliptic estimates}\label{EllipticEstimates}
Our object here is to prove Lemmas \ref{A0Lemma},
\ref{SpaceDerivativesEllipticLemma} and \ref{SpaceTimeEllipticLemma}.
\subsection{Basic estimates}
We first prove existence and uniqueness for the equation
\begin{equation}\label{EllipticEquation}
   \Delta u - \abs{\phi}^2 u = - \Im (\phi f)
\end{equation}
on $\R^4$.
\begin{lemma}\label{EllipticLemma2}
Let $\phi \in \dot H^1$ and $f \in L^2$. Then the equation
\eqref{EllipticEquation}
has a unique (real-valued) solution $u \in \dot H^1$, and
\begin{equation}\label{EllipticEstimate}
   \Sobdotnorm{u}{1} \le 2 \twonorm{f}{}.
\end{equation}
\end{lemma}
\begin{proof} Recall that $\dot H^1$, as defined in section \ref{Notation},
is a Hilbert space with inner product $\int \nabla u \cdot \overline{\nabla v}$
(by Plancherel's theorem), and that $\dot H^1 \hookrightarrow L^4$.
We denote by $\Re \dot H^1$ the corresponding
real Hilbert space, with inner product $\int \nabla u \cdot \nabla v$.

By definition, $u \in \dot H^1$ solves \eqref{EllipticEquation} in
the sense of distributions iff
\begin{equation}\label{EllipticA}
   \int_{\R^4} \bigl( \nabla u \cdot \nabla v
   + \abs{\phi}^2 u v \bigr) \, dx = \int \Im (\phi f) v \, dx
\end{equation}
for all $v \in \Schwartz$. Since $\Schwartz$ is dense in $\dot H^1$ and
\begin{align}
   \label{EllipticB}
   \abs{ \int \abs{\phi}^2 u v \, dx} &\le \Lpnorm{\phi}{4}^2 
\Lpnorm{u}{4} \Lpnorm{v}{4}
   \lesssim \Sobdotnorm{\phi}{1}^2 \Sobdotnorm{u}{1} \Sobdotnorm{v}{1},
   \\
   \label{EllipticC}
   \abs{ \int \Im (\phi f) v \, dx} &\le \Lpnorm{\phi}{4} 
\Lpnorm{f}{2} \Lpnorm{v}{4}
   \lesssim \Sobdotnorm{\phi}{1} \Sobdotnorm{v}{1} \Lpnorm{f}{2},
\end{align}
we conclude that $u$ solves \eqref{EllipticEquation} iff
\eqref{EllipticA} holds for all $v \in \dot H^1$.
Taking $v = \overline u$ gives
$$
   \Sobdotnorm{u}{1}^2 + \twonorm{\phi u}{}^2 = \int \Im (\phi f) 
\overline u \, dx
   \le \twonorm{\phi u}{} \twonorm{f}{},
$$
and since $(a+b)^2 \le 2(a^2 + b^2)$ for all $a,b \in \R$,
we conclude that $N^2 \le 2 N \twonorm{f}{}$ where
$N = \Sobdotnorm{u}{1} + \twonorm{\phi u}{} < \infty$. Therefore
\eqref{EllipticEstimate} holds, and uniqueness follows.
Of course, $u$ must be real, since if $u$ solves \eqref{EllipticEquation},
then $\Im u$ solves the same equation with $f = 0$, and therefore
$\Im u = 0$ by what we just proved.

To prove existence, observe that the left hand side of \eqref{EllipticA}
defines an inner product on $\Re \dot H^1$, and in view of \eqref{EllipticB},
the corresponding norm is equivalent to the usual norm. Moreover,
by \eqref{EllipticC}, the right hand side of \eqref{EllipticA} is a bounded
linear functional $F(v)$ on $\Re \dot H^1$. Existence therefore
follows from the Riesz representation theorem.
\end{proof}
\begin{remark}\label{CuRemark}
As discussed in the introduction, our method can be modified to 
generalize the result of Cuccagna \cite{Cu} for MKG on $\R^{1+3}$ to 
large data in $H^{s}$, $s > 3/4$. For this, we need the fact that 
\eqref{EllipticEquation} has a unique solution in $\dot H^{1}(\R^{3})$ 
for $\phi \in H^{3/4}(\R^{3})$ and $f \in H^{-1/4}(\R^{3})$.
Again we multiply the equation by $\overline u$ and integrate.
Using Plancherel's theorem we get
$$
   \Sobdotnorm{u}{1}^2 + \twonorm{\phi u}{}^2
   \le \Sobnorm{\phi u}{1/4} \Sobnorm{f}{-1/4}
   + \Sobnorm{\phi \overline u}{1/4} \Sobnorm{f}{-1/4}
$$
and since
$$
  \Sobnorm{\phi u}{1/4} \lesssim \Sobnorm{\phi}{3/4} 
  \Sobdotnorm{u}{1},
$$
on $\R^{3}$, we get $\Sobdotnorm{u}{1} \lesssim \Sobnorm{\phi}{3/4} \Sobnorm{f}{-1/4}$.
It is also easy to show that the operator $B_{0}$ defined by 
\eqref{B0Def} is bounded in $L^{2}$ for $\phi, A_{j} \in 
H^{3/4}(\R^{3})$.
\end{remark}
\medskip\medskip
Next, we prove a difference estimate for \eqref{EllipticEquation}.
\begin{lemma}\label{EllipticLemma3}
Let $\phi, \psi \in \dot H^1$ and $f, g \in L^2$. Let $u,v \in \dot H^1$
be the solutions of
\begin{align*}
   \Delta u - \abs{\phi}^2 u &= - \Im (\phi f),
   \\
   \Delta v - \abs{\psi}^2 v &= - \Im (\psi g).
\end{align*}
Then
$$
   \Sobdotnorm{u-v}{1} \lesssim
   \Sobdotnorm{\phi-\psi}{1} + \twonorm{f-g}{}
$$
where the suppressed constant is a polynomial in
$\Sobdotnorm{\phi}{1}$, $\Sobdotnorm{\psi}{1}$ and $\twonorm{g}{}$.
\end{lemma}
\begin{proof} Subtracting the equations gives
$$
   \Delta(u-v) - \abs{\phi}^2 (u-v)
   = \bigl( \abs{\phi}^2 - \abs{\psi}^2 \bigr) v
   - \Im \left[ \phi(f-g) \right]
   - \Im \left[ (\phi-\psi)g \right].
$$
Then by a density argument as in the previous proof,
\begin{align*}
   &\int \left( \nabla (u-v) \cdot \nabla (u-v) + \abs{\phi}^2 (u-v)^2 \right)
   \, dx
   \\
   &\qquad \quad
   = \int \left( \bigl( \abs{\psi}^2 - \abs{\phi}^2 \bigr) v
   + \Im \left[ \phi(f-g) \right]
   + \Im \left[ (\phi-\psi)g \right] \right) (u-v) \, dx
   \\
   &\qquad \quad
   \le
   \Lpnorm{\phi-\psi}{4} \bigl( \Lpnorm{\phi}{4} + \Lpnorm{\psi}{4} \bigr)
   \Lpnorm{v}{4} \Lpnorm{u-v}{4}
   \\
   &\qquad \qquad \qquad \quad
   + \Lpnorm{\phi}{4} \twonorm{f-g}{} \Lpnorm{u-v}{4}
   + \Lpnorm{\phi-\psi}{4} \twonorm{g}{} \Lpnorm{u-v}{4},
\end{align*}
giving the desired conclusion.
\end{proof}
We now consider the more general equation
\begin{equation}\label{GeneralEllipticEquation}
   \Delta u - \abs{\phi}^2 u = f
\end{equation}
\begin{lemma}\label{EllipticLemma} Given $\phi \in \dot H^1$ and $f 
\in L^{4/3}$, the equation \eqref{GeneralEllipticEquation}
has a unique solution $u \in \dot H^1$, and
\begin{equation}\label{EllipticSectionA}
   \Sobdotnorm{u}{1} \le C \Lpnorm{f}{4/3}
\end{equation}
where $C$ is independent of $\phi, f$ and $u$. Moreover, if
\begin{align}
   \label{EllipticSectionB}
   \Delta u - \abs{\phi}^2 u &= f,
   \\
   \label{EllipticSectionC}
   \Delta v - \abs{\psi}^2 v &= g,
\end{align}
where $u,v,\phi,\psi \in \dot H^1$ and $f,g \in L^{4/3}$, then
$$
   \Sobdotnorm{u-v}{1} \le C \bigl( \Sobdotnorm{\phi}{1} + \Sobdotnorm{\psi}{1}
   \bigr) \Lpnorm{g}{4/3} \Sobdotnorm{\phi-\psi}{1}
   + C \Lpnorm{f-g}{4/3},
$$
with the same constant $C$ as above.
\end{lemma}
\begin{proof}
Proceed as in the proof of Lemma \ref{EllipticLemma2},
but with the right hand side of \eqref{EllipticA}
replaced by $- \int vf \, dx$. Thus \eqref{EllipticC}
is replaced by
$$
   \abs{ \int v f \, dx} \le \Lpnorm{v}{4} \Lpnorm{f}{4/3}
   \lesssim \Sobdotnorm{v}{1} \Lpnorm{f}{4/3}.
$$
Existence then follows, and any $\dot H^1$ solution satisfies
$$
   \Sobdotnorm{u}{1}^2 + \twonorm{\phi u}{}^2 = - \int u f \, dx
   \le C \Sobdotnorm{u}{1} \Lpnorm{f}{4/3},
$$
where $C$ is independent of $u, f$ and $\phi$, and
\eqref{EllipticSectionA} follows.

Subtracting \eqref{EllipticSectionB} from \eqref{EllipticSectionA} gives
$$
   \bigl(\Delta - \abs{\phi}^2\bigr) (u-v) =
   \bigl( \abs{\phi}^2 - \abs{\psi}^2 \bigr) v + f-g,
$$
and applying \eqref{EllipticSectionA} gives the desired estimate.
\end{proof}
Next we prove a uniqueness result in space-time:
\begin{lemma}\label{SpaceTimeUniqueness}
Suppose
$$
  \phi \in C\bigl( [0,T] , \dot H^1\bigr)
  \quad \text{and} \quad
  u \in L^2\bigl( [0,T] , \dot H^1\bigr),
$$
and that $u$ solves
$$
  \Delta u - \abs{\phi}^2 u = 0 \quad \text{on} \quad (0,T) \times \R^4
$$
in the sense of distributions. Then $u = 0$.
\end{lemma}
\begin{proof}
Set $S_T = (0,T) \times \R^4$.
For every test function $v(t,x)$ in $C_c^\infty(S_T)$,
\begin{equation}\label{WeakSolution}
  \int \bigl\{ \nabla u \cdot \nabla v + \abs{\phi}^2 uv \bigr\}
  \, dt \, dx = 0.
\end{equation}
The left hand side is a bounded linear functional in $v$. In fact,
\begin{align*}
  \abs{ \int \nabla u \cdot \nabla v \, dt \, dx }
  &\le \twonorm{\nabla u}{(S_T)} \twonorm{\nabla v}{(S_T)}
  = \norm{u}_{L_t^2(\dot H^1)} \norm{v}_{L_t^2(\dot H^1)},
  \\
  \abs{\int \abs{\phi}^2 uv \, dt \, dx }
  &\lesssim
  \norm{\phi}_{L_t^\infty(\dot H^1)}^2
  \norm{u}_{L_t^2(\dot H^1)}
  \norm{v}_{L_t^2(\dot H^1)}.
\end{align*}
Here we used H\"older's inequality and the embedding
$\dot H^1 \hookrightarrow L^4$.

But $C_c^\infty(S_T)$ is dense in $L^2( [0,T] , \dot H^1 )$,
so it follows that \eqref{WeakSolution} must hold for all
$v \in L^2( [0,T] , \dot H^1 )$. Taking $v = \overline u$ gives
$$
 \int \bigl\{ \abs{\nabla u}^2 + \abs{\phi}^2 \abs{u}^2 \bigr\}
  \, dt \, dx = 0.
$$
This implies $\nabla u = 0$, hence $u = 0$ ($\dot H^1$, as we have
defined it, does not contain any nonzero constants).
\end{proof}
\subsection{Higher regularity estimates}
Suppose
$$
  \phi \in C\bigl( [0,T] , \dot H^1 \bigr)
  \quad \text{and} \quad
  f \in C\bigl( [0,T] , L^{4/3} \bigr).
$$
By Lemma \ref{EllipticLemma}, the equation
\begin{equation}\label{EllipticEq}
  \Delta u - \abs{\phi}^2 u = f
\end{equation}
has a unique solution
$$
  u \in C\bigl( [0,T] , \dot H^1\bigr).
$$
We shall prove the following higher regularity estimate.
\begin{lemma}\label{SpaceTimeHigherOrderEllipticLemma}
Let $m,M$ be non-negative integers. If
$$
  \partial_t^j \partial_x^\alpha \phi \in C\bigl( [0,T] , \dot H^1 \bigr)
  \quad \text{and} \quad
  \partial_t^j \partial_x^\alpha f \in C\bigl( [0,T] , L^{4/3} \bigr)
$$
for all $j \le m$ and $\abs{\alpha} \le M$, then
$$
  \partial_t^j \partial_x^\alpha u \in C\bigl( [0,T] , \dot H^1 \bigr)
$$
for $j \le m$ and $\abs{\alpha} \le M$, and
\begin{multline*}
  \bignorm{\partial_t^j \partial_x^\alpha u}_{L_t^\infty(\dot H^1)}
  \lesssim
  \bignorm{\partial_t^j \partial_x^\alpha \phi}_{L_t^\infty(\dot H^1)}
  \norm{\phi}_{L_t^\infty(\dot H^1)}
  \norm{f}_{L_t^\infty(L^{4/3})}
  \\
  +
  \bignorm{\partial_t^j \partial_x^\alpha f}_{L_t^\infty(L^{4/3})}
  +
  \eta_{\alpha},
\end{multline*}
where $\eta_{\alpha}$ is a lower order term which depends
continuously on the norms
$$
  \bignorm{\partial_t^k \partial_x^\beta \phi}_{L_t^\infty(\dot H^1)}
  \quad \text{and} \quad
  \bignorm{\partial_t^k \partial_x^\beta f}_{L_t^\infty(L^{4/3})}
$$
for all $k \le j$ and $\abs{\beta} \le \abs{\alpha}$
satisfying $k + \abs{\beta} < j + \abs{\alpha}$.
Here all $L^\infty$-norms are taken over $[0,T]$.
\end{lemma}
The proof is by induction on $m$ and $M$.
Denote by $P(m,M)$ the statement that the lemma
holds for the pair $(m,M)$. Since $P(0,0)$ is true by Lemma \ref{EllipticLemma},
it is enough, by an obvious induction, to prove
\begin{subequations}\label{Inductive}
\begin{gather}
  \label{Inductive1}
  P(m,0) \implies P(m+1,0),
  \\
  \label{Inductive2}
  P(0,M) \implies P(0,M+1),
  \\
  \label{Inductive3}
  P(m+1,M), P(m,M+1) \implies P(m+1,M+1).
\end{gather}
\end{subequations}
The key to proving these implications is the following:
\begin{lemma}\label{SpaceTimeOneDerivativeEllipticLemma}
If, for some $0 \le \mu \le 4$,
$$
  \partial_\mu \phi \in C\bigl( [0,T] , \dot H^1 \bigr)
  \quad \text{and} \quad
  \partial_\mu f \in C\bigl( [0,T] , L^{4/3} \bigr),
$$
then $\partial_\mu u \in C( [0,T] , \dot H^1 )$ and
$$
  \norm{\partial_\mu u}_{L_t^\infty(\dot H^1)}
  \lesssim
  \norm{f}_{L_t^\infty(L^{4/3})}
  \norm{\phi}_{L_t^\infty(\dot H^1)}
  \norm{\partial_\mu \phi}_{L_t^\infty(\dot H^1)}
  \\
  + \norm{\partial_\mu f}_{L_t^\infty(L^{4/3})},
$$
where the $L^\infty$-norms are taken over $[0,T]$.
\end{lemma}
Before proving Lemma \ref{SpaceTimeOneDerivativeEllipticLemma},
let us use it to prove \eqref{Inductive}.
\paragraph{Proof of \eqref{Inductive1}.}
Apply $\partial_t^m$ to both sides of \eqref{EllipticEq}. This gives
\begin{equation}\label{Bigsum1}
   \bigl(\Delta - \abs{\phi}^2 \bigr) \partial_t^m u
   = \sum_{\text{$j+k+l = m$, $l < m$}}
   c_{j k l}
   (\partial_t^j \phi)
   ( \overline{ \partial_t^k \phi })
   (\partial_t^l u) + \partial_t^m f.
\end{equation}
Denote the right hand side of this equation by $F$.
By Lemma \ref{SpaceTimeOneDerivativeEllipticLemma},
if we can show that $F$ and $\partial_t F$
belong to $C([0,T],L^{4/3})$, then it follows that
$$
  \partial_t^{m+1} u \in C([0,T],\dot H^1)
$$
and
$$
  \norm{\partial_t^{m+1} u}_{L_t^\infty(\dot H^1)}
  \lesssim
  \norm{F}_{L_t^\infty(L^{4/3})}
  \norm{\phi}_{L_t^\infty(\dot H^1)}
  \norm{\partial_t \phi}_{L_t^\infty(\dot H^1)}
  \\
  + \norm{\partial_t F}_{L_t^\infty(L^{4/3})}.
$$
But by H\"older's inequality and
the embedding $\dot H^1 \hookrightarrow L^4$,
\begin{multline*}
  \norm{F}_{L_t^\infty(L^{4/3})}
  \lesssim
  \sum
  \bignorm{\partial_t^j \phi}_{L_t^\infty(\dot H^1)}
  \norm{\partial_t^k \phi}_{L_t^\infty(\dot H^1)}
  \norm{\partial_t^l u}_{L_t^\infty(\dot H^1)}
  \\
  + \bignorm{\partial_t^m f}_{L_t^\infty(L^{4/3})},
\end{multline*}
and using the hypothesis $P(m,0)$ to bound
$\norm{\partial_t^l u}_{L_t^\infty(\dot H^1)}$,
we conclude that $\norm{F}_{L_t^\infty(L^{4/3})}$
is bounded by a continuous function of the norms
$\bignorm{\partial_t^j \phi}_{L_t^\infty(\dot H^1)}$
and $\bignorm{\partial_t^j f}_{L_t^\infty(L^{4/3})}$
for all $j \le m$.

Next apply $\partial_t$ to $F$.
Since $l < m$ in \eqref{Bigsum1},
we can again use the hypothesis $P(m,0)$ to bound
$\norm{\partial_t^{l+1} u}_{L_t^\infty(\dot H^1)}$.
Thus, arguing as before,
\begin{multline*}
  \norm{\partial_t F}_{L_t^\infty(L^{4/3})}
  \lesssim
  \norm{\partial_t^{m+1} \phi}_{L_t^\infty(\dot H^1)}
  \norm{\phi}_{L_t^\infty(\dot H^1)}
  \norm{f}_{L_t^\infty(L^{4/3})}
  \\
  + \norm{\partial_t^{m+1} f}_{L_t^\infty(L^{4/3})}
  + \text{l.o.t.},
\end{multline*}
where l.o.t. stands for a lower order term
depending continuously on the norms
$\bignorm{\partial_t^j \phi}_{L_t^\infty(\dot H^1)}$
and $\bignorm{\partial_t^j f}_{L_t^\infty(L^{4/3})}$
for all $j \le m$.

Thus $P(m+1,0)$ holds, completing the proof of \eqref{Inductive1}.
The proof of \eqref{Inductive2} is quite similar and is omitted.
\paragraph{Proof of \eqref{Inductive3}.}
Let $\abs{\alpha} = M+1$.
Applying $\partial_t^m \partial_x^\alpha$
to both sides of \eqref{EllipticEq} gives
\begin{equation}\label{Bigsum}
   \bigl(\Delta - \abs{\phi}^2 \bigr) \partial_t^m \partial_x^{\alpha} u
   = \sum
   c_{\beta \gamma \delta j k l}
   (\partial_t^j \partial_x^\beta \phi)
   ( \overline{ \partial_t^k \partial_x^\gamma \phi })
   (\partial_t^l \partial_x^\delta u) + \partial_t^m \partial_x^\alpha f,
\end{equation}
where the sum is over all non-negative integers $j,k,l$ and multi-indices
$\beta,\gamma,\delta$ such that
$$
  \beta+\gamma+\delta = \alpha,
  \quad j + k + l = m, \quad l + \abs{\delta} \le m + M.
$$
Let $F$ be the right hand side of \eqref{Bigsum}.
If $F$ and $\partial_t F$
belong to $C([0,T],L^{4/3})$,
then by Lemma \ref{SpaceTimeOneDerivativeEllipticLemma},
$$
  \partial_t^{m+1} \partial_x^\alpha u \in C([0,T],\dot H^1)
$$
and
$$
  \norm{\partial_t^{m+1} \partial_x^\alpha u}_{L_t^\infty(\dot H^1)}
  \lesssim
  \norm{F}_{L_t^\infty(L^{4/3})}
  \norm{\phi}_{L_t^\infty(\dot H^1)}
  \norm{\partial_t \phi}_{L_t^\infty(\dot H^1)}
  \\
  + \norm{\partial_t F}_{L_t^\infty(L^{4/3})}.
$$
Since $F$ and $\partial_t F$ are estimated just as in the
proof of \eqref{Inductive1}, we will not go into details.
The key point is that since $l \le m$, $\abs{\delta} \le M+1$
and $l + \abs{\delta} \le m + M$
in \eqref{Bigsum}, the hypotheses $P(m+1,M)$ and $P(m,M+1)$
allow us to bound
$\norm{\partial_t^l \partial_x^\delta u}_{L_t^\infty(\dot H^1)}$
and
$\norm{\partial_t^{l+1} \partial_x^\delta u}_{L_t^\infty(\dot H^1)}$.
\paragraph{Proof of Lemma \ref{SpaceTimeOneDerivativeEllipticLemma}.}
Under the hypotheses of the lemma,
\begin{equation}\label{L2inTime}
  \partial_\mu u \in L^2\bigl([0,T],\dot H^1\bigr).
\end{equation}
Before proving this, let us show that it implies the conclusion
of the lemma.

Indeed, $u$ solves \eqref{EllipticEq} in the sense
of distributions on $S_T = (0,T) \times \R^4$, and if we apply
$\partial_\mu$ to both sides, it follows that
\begin{equation}\label{SpaceTimeEquation1}
  \bigl( \Delta - \abs{\phi}^2 \bigr) \partial_\mu u
  = \partial_\mu\bigl(\abs{\phi}^2\bigr) u + \partial_\mu f.
\end{equation}
in the sense of distributions on $S_T$. (The use of the product
rule for derivatives is easily justified in view of \eqref{L2inTime}.)
Denote by $F$ the right hand side of the last equation. Then
$$
  F = 2 \Re \bigl( \phi \overline{\partial_\mu \phi} \bigr) u
  + \partial_\mu f,
$$
and so
$$
  \norm{F}_{L_t^\infty(L^{4/3})}
  \lesssim
  \norm{\phi}_{L_t^\infty(\dot H^1)}
  \norm{\partial_\mu \phi}_{L_t^\infty(\dot H^1)}
  \norm{u}_{L_t^\infty(\dot H^1)}
  + \norm{\partial_\mu f}_{L_t^\infty(L^{4/3})}.
$$
It then follows by Lemma \ref{EllipticLemma} that
the equation
\begin{equation}\label{SpaceTimeEquation2}
  \bigl( \Delta - \abs{\phi}^2 \bigr) v = F
\end{equation}
has a solution $v \in C( [0,T] , \dot H^1 )$, and
$$
  \norm{v}_{L_t^\infty(\dot H^1)}
  \lesssim
  \norm{f}_{L_t^\infty(L^{4/3})}
  \norm{\phi}_{L_t^\infty(\dot H^1)}
  \norm{\partial_\mu \phi}_{L_t^\infty(\dot H^1)}
  \\
  + \norm{\partial_\mu f}_{L_t^\infty(L^{4/3})}.
$$
Subtracting \eqref{SpaceTimeEquation2} from
\eqref{SpaceTimeEquation1} gives
$$
  \bigl( \Delta - \abs{\phi}^2 \bigr) (\partial_\mu u - v) = 0.
$$
Thus $\partial_\mu u = v$ by Lemma \ref{SpaceTimeUniqueness},
proving the conclusion of Lemma \ref{SpaceTimeOneDerivativeEllipticLemma}.

It remains to prove \eqref{L2inTime}.
For technical reasons, we fix $0 < t_0 < T$ and prove
\eqref{L2inTime} with the interval $[0,T]$ replaced by $[0,t_0]$.
A similar argument works for the interval $[t_0,T]$,
giving the statement in the entire interval $[0,T]$.

We shall require the following facts about the difference quotients
$$
   \Delta_h^j u(t,x) = \frac{ u(t,x+he_j) - u(t,x) }{h}
   \quad \text{and} \quad
   \Delta_h^0 u(t,x) = \frac{ u(t+h,x) - u(t,x) }{h},
$$
where $e_1,\dots,e_4$ are the standard basis vectors of $\R^4$.
\begin{lemma}\label{DiffLemmaA}
If $f$ belongs to $C( [0,T] , L^{4/3})$
and the distributional derivative
$\partial_\mu f$ belongs to $C( [0,T] , L^{4/3} )$,
for some $0 \le \mu \le 4$, then
$$
  \norm{\Delta_h^\mu f}_{L^\infty([0,t_0], L^{4/3})}
  = O(1)
  \quad \text{as} \quad h \to 0^+.
$$
Moreover, the same conclusion holds with $f$ replaced by
$\phi$ and $L^{4/3}$ by $\dot H^1$.
\end{lemma}
\begin{proof} We have
$$
  \Delta_h^0 f(t) = \frac{1}{h} \int_0^{h}
  \partial_t f(t+s) \, ds \qquad \text{($L^{4/3}$-valued integral)}
$$
whence
$$
  \Lpnorm{\Delta_h^0 f(t)}{4/3}
  \le \sup_{0 \le s \le h} \Lpnorm{\partial_t f(t+s)}{4/3}.
$$
The same proof works for $\dot H^1$.
If $1 \le j \le 4$, then (discarding the time variable)
$$
  \Delta_h^j f(x) = \frac{1}{h} \int_0^{h}
  \partial_j f(x+se_j) \, ds.
$$
Thus, by Minkowski's integral inequality and
the translation invariance of the norm,
$$
  \bigLpnorm{\Delta_h^j f}{4/3}
  \le \Lpnorm{\partial_j f}{4/3}.
$$
This is certainly valid for smooth $f$, and hence in general
by using an approximation of the identity. For $\dot H^1$
we write
$$
  \bigSobdotnorm{\Delta_h^j \phi}{1}^2
  = \int \abs{\xi}^2 \abs{\frac{e^{ih\xi_j}-1}{h}}^2
  \bigabs{\widehat \phi}^2 \, d\xi
$$
and note that $\abs{e^{ih\xi_j}-1} \le h \abs{\xi_j}$.
It follows that $\bigSobdotnorm{\Delta_h^j \phi}{1}
\le \Sobdotnorm{\partial_j \phi}{1}$.
\end{proof}
We are now ready to prove \eqref{L2inTime}.
By the difference estimate in Lemma \ref{EllipticLemma},
\begin{multline*}
  \norm{\Delta_h^\mu u}_{L^\infty([0,t_0],\dot H^1)}
  \lesssim
  \norm{f}_{L^\infty([0,T],L^{4/3})}
  \norm{\phi}_{L^\infty([0,T],\dot H^1)}
  \norm{\Delta_h^\mu \phi}_{L^\infty([0,t_0],\dot H^1)}
  \\
  + \norm{\Delta_h^\mu f}_{L^\infty([0,t_0],L^{4/3})},
\end{multline*}
for all $0 < h < T-t_0$.
In view of Lemma \ref{DiffLemmaA}, the right hand
side is $O(1)$ as $h \to 0^+$. Applying H\"older's inequality in time
then gives
$$
  \norm{\Delta_h^\mu u}_{L^2([0,t_0],\dot H^1)}
  = O(1) \quad \text{as} \quad h \to 0^+,
$$
so by weak compactness, there is a sequence $h_j \to 0$
such that $\Delta_{h_j}^\mu u$ converges weakly in
$L^2([0,t_0],\dot H^1)$ to some limit $v$ as $j \to \infty$.
But this implies that $\Delta_{h_j}^\mu u \to v$ also in
the sense of distributions on $(0,t_0) \times \R^4$.
On the other hand, we know that
$\Delta_h^\mu u \to \partial_\mu u$
in the distributional sense as $h \to 0$, and so we conclude that
$\partial_\mu u = v$. This proves \eqref{L2inTime}
on the interval $[0,t_0]$.
\subsection{Proofs of  Lemmas \ref{A0Lemma},
\ref{SpaceDerivativesEllipticLemma} and \ref{SpaceTimeEllipticLemma}}
\label{EllipticLemmaProofs}
First, Lemma \ref{A0Lemma} is an immediate
corollary of Lemmas \ref{EllipticLemma2} and \ref{EllipticLemma3}.

Secondly, to prove Lemma \ref{SpaceDerivativesEllipticLemma}, we apply
Lemma \ref{SpaceTimeHigherOrderEllipticLemma} with $m = 0$,
$M = k$ and
$$
  f = - \Im \bigl(\phi \overline{\partial_t \phi} \bigr).
$$
Since $\Lambda^k \phi \in \X_2 = \scrH^{s,\theta}$ by hypothesis, we have,
in view of \eqref{EnergyEmbedding},
\begin{equation}\label{PhiReg}
  \partial_x^\alpha \phi \in C\bigl( [0,T] , \dot H^1 \bigr)
  \quad \text{and} \quad
  \partial_t \partial_x^\alpha \phi \in C\bigl( [0,T] , L^2 \bigr)
  \quad \text{for} \quad \abs{\alpha} \le k.
\end{equation}
Thus, it suffices to check that
$\partial_x^\alpha f \in C( [0,T] , L^{4/3} )$
for $\abs{\alpha} \le k$, and
\begin{equation}\label{fEstimate}
  \norm{\partial_x^\alpha f}_{L_t^\infty(L^{4/3})}
  \lesssim \norm{\phi}_{\X_2} \norm{\partial_x^\alpha \phi}_{\X_2}
  + \eta\Bigl( \sum_{\abs{\beta} < \abs{\alpha}}
  \norm{\partial_x^\beta \phi}_{\X_2} \Bigr),
\end{equation}
where $\eta$ is continouous. If $\alpha = 0$, we have,
using H\"older's inequality and the embedding
$\dot H^1 \hookrightarrow L^4$,
$$
  \norm{f}_{L^{4/3}}
  \le \norm{\phi}_{L^4} \norm{\partial_t \phi}_{L^2}
  \lesssim \norm{\phi}_{\dot H^1} \norm{\partial_t \phi}_{L^2}
$$
uniformly in $0 \le t \le T$, and by \eqref{EnergyEmbedding},
$$
  \norm{\phi}_{\dot H^1}, \norm{\partial_t \phi}_{L^2}
  \lesssim \norm{\phi}_{\X_2},
$$
giving \eqref{fEstimate} for $\alpha = 0$. When $\alpha \neq 0$
one can apply the product rule and estimate each term
as above. We leave the details to the interested reader.

Finally, Lemma \ref{SpaceTimeEllipticLemma} is also
proved by an application of Lemma
\ref{SpaceTimeHigherOrderEllipticLemma}.
We are given non-negative integers $m,M$ such that
$$
  \partial_t^j \partial_x^\alpha \phi \in  C\bigl([0,T],L^2\bigr)
  \quad \text{for all} \quad j \le m+1 \quad \text{and all} \quad 
  \abs{\alpha} \le M+1.
$$
Again we set
$$
  f = - \Im \bigl(\phi \overline{\partial_t \phi} \bigr).
$$
By Lemma \ref{SpaceTimeHigherOrderEllipticLemma}
it suffices to check that
$$
  \partial_t^j \partial_x^\alpha f \in  C\bigl([0,T],L^{4/3}\bigr)
  \quad \text{for all} \quad j \le m \quad \text{and all} \quad 
  \abs{\alpha} \le M.
$$
Again, one simply applies the product rule for derivatives and
estimates each term as in the proof of \eqref{fEstimate}.
\section*{Appendix}
Here we prove \eqref{Tataru}. First,
$$
   \bigmixednorm{\Lambda^{s-1}(-\Delta)^{-1}(uv)}{p}{q}
   \lesssim \bigmixednorm{(-\Delta)^{-1}(uv)}{p}{q}
   + \bigmixednorm{(-\Delta)^{\frac{s-3}{2}}(uv)}{p}{q}
$$
by Lemma \ref{BesselRieszLemma}, so it suffices to show
that the two terms on the right hand side are both
$\lesssim \spacetimenorm{u}{s}{\theta} \spacetimenorm{v}{s-1}{\theta}$.
To this end, we apply the following theorem (stated here for $\R^{1+4}$ only)
of Klainerman-Tataru \cite{Kl-Ta}:
\begin{plaintheorem} Let $1 \le p \le \infty$, $1 \le q < \infty$ and
set $\gamma = 2  - \frac{1}{2p} - \frac{2}{q}$. Assume that
\begin{gather}
   \label{A1}
   \frac{1}{p} \le \frac{3}{2}\left(1 - \frac{1}{q}\right),
   \\
   \label{A2}
   0 < \sigma < 4  - \frac{2}{p} - \frac{4}{q},
   \\
   \label{A3}
   s_1,s_2 < \gamma,
   \\
   \label{A4}
   s_1 + s_2 + \sigma = 2\gamma.
\end{gather}
Then
$$
   \mixednorm{(-\Delta)^{-\sigma/2} (uv)}{p}{q} \lesssim
   \spacetimenorm{u}{s_1}{\theta} \spacetimenorm{v}{s_2}{\theta},
$$
where $\theta > \half$.
\end{plaintheorem}
\begin{plainremarks}
(1) When $s_1 = s_2$ this follows from Theorem 4 in \cite{Kl-Ta}
(see also \cite[Principle 3.2]{Kl-Se}).
In \cite{Kl-Ta}, however, the estimate was stated using the
space-time fractional derivative operator
$(-\Delta_{t,x})^{-\sigma/2}$. Nevertheless, an inspection of
their proof shows that it works equally well for
$(-\Delta)^{-\sigma/2}$ (see \cite[Chapter 2]{Se2}).
In our statement of the theorem we have also included
the end-point case due to Keel and Tao \cite{Ke-Ta}, although
we do not use this.

\medskip
\noindent
(2) The asymmetric case $s_1 \neq s_2$ is derived
as in the proof of Theorem 5 in \cite{Kl-Ta}.
(The statement of that theorem contains the
condition (in our notation) $\sigma \le \gamma$, but an inspection of
the proof shows that this is superfluous.)
Let us just give a heuristic explanation of why
the asymmetric case essentially reduces to
the symmetric situation. Rewrite the estimate as follows:
$$
  \mixednorm{D^{-\sigma} (D^{-s_1}u \cdot D^{-s_2} v)}{p}{q} \lesssim
  \spacetimenorm{u}{0}{\theta} \spacetimenorm{v}{0}{\theta},
$$
where $D^\alpha = (-\Delta)^{\alpha/2}$.
Denote by $\xi$ and $\eta$ the Fourier frequencies of $u$ and $v$
corresponding to the spatial variable $x$. Then the frequency of
the product $uv$ is $\xi + \eta$, and in Fourier space,
$D^{-\sigma} (D^{-s_1}u \cdot D^{-s_2} v)$ is a weighted convolution,
with weights
$$
  \frac{1}{\abs{\xi+\eta}^\sigma \abs{\xi}^{s_1} \abs{\eta}^{s_2}}
$$
The idea is that the weights can be redistributed so as to
get equal powers of $\abs{\xi}$ and $\abs{\eta}$.
This is obviously possible if the frequencies of $u$ and
$v$ are comparable.
If, on the other hand, $\abs{\xi} \gg \abs{\eta}$, say,
then $\abs{\xi+\eta} \sim \abs{\xi}$, and so
$$
  \frac{1}{\abs{\xi+\eta}^\sigma \abs{\xi}^{s_1} \abs{\eta}^{s_2}}
  \lesssim
  \frac{1}{\abs{\xi}^{\gamma} \abs{\eta}^{\gamma}}
$$
provided $s_1,s_2 \le \gamma$ (recall that $s_1 + s_2 + \sigma = 2\gamma$).
Thus we are in the case $\sigma = 0$ and $s_1 = s_2 = \gamma$,
which by H\"older's inequality is reduced to a linear
Strichartz estimate.
\end{plainremarks}
Now let $p$ and $q$ be defined as in section \ref{MainEstimates}.
Using the definition of $p$ in \eqref{prdef}, we
see that \eqref{A1} is equivalent to
$1/q \le (2/3)(\theta + 2\varepsilon)$, and the latter
evidently holds, since $1/q \le 1/4$
by \eqref{qdef}, \eqref{Parameters1} and the assumption $s < 2$.
Thus \eqref{A1} holds.

Next we have to check that \eqref{A2} holds with
$\sigma = 2$ (then it also holds with $\sigma = 3-s$, of course),
but using the definition of $p$ in \eqref{prdef}, we find that
\eqref{A2} is equivalent to
$$
  \frac{4}{q} < 2\theta - 1 + 4 \varepsilon,
$$
which is true by \eqref{qdef}.

Now set
$s_1 = 1 + \delta$ and $s_2 = \delta$,
where we have defined
\begin{equation}\label{A5}
   \delta = \frac{3-\sigma}{2} - \frac{1}{2p} - \frac{2}{q}.
\end{equation}
With this choice, \eqref{A4} clearly holds. Note that \eqref{A3} holds provided
$$
   \delta < 1 - \frac{1}{2p} - \frac{2}{q}.
$$
It suffices to check this when $\sigma = 3 - s$ (then it also
holds for $\sigma = 2$), but in this case it is obvious
since $s < 2$ and, from \eqref{A5},
\begin{equation}\label{delta}
   \delta = \frac{s}{2} - \frac{1}{2p} - \frac{2}{q}.
\end{equation}

It remains to check that $s_1 \le s$ and $s_2 \le s-1$.
This is equivalent to $\delta \le s-1$,
and again we only have to check this for $\sigma = 3 - s$,
in which case it reduces to, by \eqref{delta},
$$
  \frac{s}{2} - \frac{1}{2p} - \frac{2}{q} \le s-1.
$$
In fact,
$$
   \frac{s}{2} - \frac{1}{2p} < s-1,
$$
for by \eqref{prdef}, this is equivalent to
$1/4 + \theta/2 + \varepsilon < s/2$, which holds by \eqref{Parameters2}.

\begin{thebibliography}{1}
%
\bibitem{Cu}
S. Cuccagna,
\textit{On the local existence for the Maxwell-Klein-Gordon system in
${R}\sp {3+1}$}, Comm. PDE \textbf{24} (1999),
no. 5-6, 851--867
%
\bibitem{E-M}
D. Eardley and V. Moncrief,
\textit{The global existence of Yang-Mills-Higgs fields in
${R}\sp {3+1}$}, Comm. Math. Phys. \textbf{83} (1982),
171--212
%
\bibitem{Ke-Ta}
M. Keel and T. Tao,
\textit{Endpoint Strichartz estimates}, Amer. J. Math. \textbf{120}
(1998), no. 5, 955--980
%
\bibitem{Kl-Ma0.2}
S. Klainerman and M. Machedon,
\textit{On the Maxwell-Klein-Gordon equation with finite energy}, Duke
Math. J. \textbf{74} (1994), 19--44
%
\bibitem{Kl-Ma2}
S. Klainerman and M. Machedon,
\textit{Remark on Strichartz type inequalities},
Int. Math. Res. Not., no. 5 (1996), 201--220
%
\bibitem{Kl-Ma5}
S. Klainerman and M. Machedon,
\textit{On the optimal local regularity for gauge field theories},
Differential and Integral Equations \textbf{10} (1997), 1019--1030
%
\bibitem{Kl-Se}
S. Klainerman and S. Selberg,
\textit{Bilinear estimates and applications to nonlinear wave equations},
to appear in Comm. in Contemporary Math.
%
\bibitem{Kl-Ta}
S. Klainerman and D. Tataru,
\textit{On the optimal local regularity for Yang-Mills equations in
$\R^{4+1}$}, J. Amer. Math. Soc. \textbf{12} (1999), 93--116
%
\bibitem{Po-Si} G. Ponce and T. Sideris, \textit{Local regularity of nonlinear
wave equations in three space dimensions}, Comm. PDE \textbf{18}
(1993), 169--177
%
\bibitem{Li1}
H. Lindblad,
\textit{Counterexamples to local existence for semilinear wave equations},
Amer. J. Math. \textbf{118} (1996), 1--16
%
\bibitem{Se2}
S. Selberg,
\textit{Multilinear space-time estimates and applications to local
existence theory for nonlinear wave equations}, Ph.D. Thesis, Princeton
University 1999
%
\bibitem{Se}
S. Selberg,
\textit{On an estimate for the wave equation and applications to
nonlinear problems}, Differential and Integral Equations
\textbf{2} (2002), 213--236
%
\bibitem{So}
C. D. Sogge, ``Lectures on nonlinear wave equations'',
Monographs in Analysis, II. International Press, 1995
%
\bibitem{St}
E. Stein, ``Singular Integrals and Differentiability Properties of Functions'',
Princeton University Press, 1970
%
\bibitem{TaoWM1}
T. Tao,
\textit{Global regularity of wave maps I. Small critical Sobolev norm 
in high dimension}, IMRN \textbf{7} (2001), 299--328
%
\bibitem{TaoWM2}
T. Tao,
\textit{Global regularity of wave maps II. Small energy in two 
dimensions}, to appear, Comm. Math. Phys.
%
\bibitem{TaoYM}
T. Tao,
\textit{Local well-posedness of the Yang-Mills equation in the 
temporal gauge below the energy norm}, preprint
\end{thebibliography}
\end{document}